\RequirePackage{fix-cm}
\documentclass[twocolumn]{svjour3}          
\smartqed  
\usepackage{graphicx}
\usepackage[square,numbers]{natbib} 
\usepackage{amsmath, amsfonts}
\usepackage{hyperref}
\renewcommand{\tens}[1]{\boldsymbol{\mathsf{#1}}}
\newcommand{\n}{\vec{n}}
\renewcommand{\d}{\,\mathrm{d}}
\newcommand{\real}{\mathbb{R}}

\newcommand{\res}{\mathcal{R}}
\newcommand{\ncell}{N}
\newcommand{\ncomp}{n_c}
\newcommand{\ndof}{n_{\mathrm{dof}}}

\newcommand{\cell}{\Omega}
\newcommand{\face}{\Gamma}
\newcommand{\neigh}{\mathcal{N}}
\newcommand{\upstr}{\mathcal{U}}
\newcommand{\downstr}{\mathcal{D}}
\newcommand{\order}{\mathcal{P}}
\newcommand{\upw}{u}

\newcommand{\acc}{\mathcal{A}}
\newcommand{\flux}{\mathcal{F}}
\newcommand{\src}{\mathcal{Q}}

\graphicspath{{./figures/}}

\usepackage{tikz}
\usetikzlibrary{calc,tikzmark,shapes,positioning}
\tikzset{every picture/.style={remember picture}}
\definecolor{re}{rgb}{0.8500, 0.3250, 0.0980}
\definecolor{bl}{rgb}{0.0000, 0.4470, 0.7410}

\tikzstyle{description} = [draw, rectangle, rounded corners, fill=white, anchor = west, text badly centered, inner sep=3pt, font = \small]

\begin{document}

\title{Efficient Reordered Nonlinear Gauss--Seidel Solvers With Higher Order For Black-Oil Models
}


\author{{\O}ystein S. Klemetsdal \and
        Atgeirr F. Rasmussen     \and
        Olav M{\o}yner           \and
        Knut-Andreas Lie
}


\institute{{\O}.S. Klemetsdal \at
            Tel.: +47 98 43 86 39\\
            Norwegian University of Science and Technology, Norway \\
            \email{oystein.klemetsdal@ntnu.no}
            \and
            A.F. Rasmussen \at
            SINTEF Digital, Norway \\
            \email{atgeirr.rasmussen@sintef.no}
            \and
            O. M{\o}yner \at
            Norwegian University of Science and Technology/SINTEF Digital, Norway \\
            \email{olav.moyner@sintef.no}
            \and
            K.-A. Lie \at
            Norwegian University of Science and Technology/SINTEF Digital, Norway \\
            \email{knut-andreas.lie@sintef.no}
}

\date{Received: date / Accepted: date}

\maketitle

\begin{abstract}
The fully implicit method is the most commonly used approach to solve black-oil problems in reservoir simulation. The method requires repeated linearization of large nonlinear systems and produces ill-condi\-tioned linear systems. We present a strategy to reduce computational time that relies on two key ideas: (\textit{i}) a sequential formulation that decouples flow and transport into separate subproblems, and (\textit{ii}) a highly efficient Gauss--Seidel solver for the transport problems. This solver uses intercell fluxes to reorder the grid cells according to their upstream neighbors, and groups cells that are mutually dependent because of counter-current flow into local clusters. The cells and local clusters can then be solved in sequence, starting from the inflow and moving gradually downstream, since each new cell or local cluster will only depend on upstream neighbors that have already been computed. Altogether, this gives optimal localization and control of the nonlinear solution process. 

This method has been successfully applied to real-field problems using the standard first-order finite volume discretization. Here, we extend the idea to first-order dG methods on fully unstructured grids. We also demonstrate proof of concept for the reordering idea by applying it to the full simulation model of the Norne oil field, using a prototype variant of the open-source OPM Flow simulator.
\end{abstract}
 
\section{Introduction}

Reservoir simulation requires the solution of large systems of nonlinear partial differential equations to compute fluid pressure, phase saturations, and component concentrations/mass fractions. In industry-grade simulations, it is most common to use implicit temporal discretization to overcome severe CFL restrictions arising because of grids with high aspect ratios, orders of magnitude variations in petrophysical properties, and large variations in fluid velocities from stagnant and slow-moving flow regions to high-flow regions around wells. For cases with weak coupling between the fluid pressure and the transport of conserved quantities, one can observe significant reduction in computational costs by splitting the flow model into a pressure equation and a set of transport equations and solve them sequentially \citep{Trangenstein1989,Watts1986}. To this end, it is common to use specialized solvers since the two subproblems typically have very different mathematical character: pressure equations are often close to elliptic, whereas transport equations have a strong hyperbolic character; see e.g., \citet{Bell1986} and \citet{Lie2016MRST-book}.

The pressure equation is usually discretized by a first-order, two-point flux-approximation scheme, giving a discrete problem that can be efficiently solved by a multigrid \citep{Trottenberg2000, Gries2014} or multiscale method \citep{Lie2017ms-review}. For the transport equations, a number of methods aim to accelerate computations by utilizing inherent locality and co-current flow properties of hyperbolic equations. Examples include streamline simulation \citep{Datta-Gupta2007, Bratvedt1996}, \textit{a priori} estimation of nonzero update regions \citep{Sheth2017}, and use of interface-localized trust regions to determine safe saturation updates \citep{Moyner2016, Klemetsdal2018}. 

Reordering methods rely on the important observation that the transport of fluid phases is always unidirectional along streamlines in the absence of gravity and capillary forces. Consequently, it is possible to order the grid cells so that the discretization matrix for the transport equations becomes lower triangular and can be solved very efficiently in a cell-by-cell manner. A branch of these methods builds upon the Cascade method \citep{Appleyard1982} and reorders grid cells based on the fluid potential \citep{Kwok2007, Shahvali2013}. Another branch uses topological traversal of the intercell flux graph \citep{Natvig2008,Lie2012,Lie2014}. The latter has the advantage that it is easy to implement using methods from standard graph theory. When gravity and capillary forces are present, the flow field will have regions of counter-current flow, which show up as cycles of mutually connected cells in the discretized flow equations. If the cells in these cycles are grouped into supernodes and solved for simultaneously, we can still use the same reordering idea.

Numerical smearing and grid-orientation effects are well-known problems in reservoir simulation. These effects are particularly evident for displacement fronts with little or no self-sharpening effects, as are often found in multicomponent models describing various mechanisms for enhanced oil recovery, e.g., polymer injection and miscible flooding. We can mitigate these effects by increasing the grid resolution and/or the formal order of the spatial discretization. Both approaches increase the computational cost by increasing the number of unknowns and/or the nonlinearity of the discretized equations, and it is therefore even more important to have a nonlinear solver with high efficacy. Most high-resolution methods---i.e., methods that deliver high formal order on smooth parts of the solution and stable propagation of discontinuities---rely on some kind of spatial reconstruction using cell-based averaged quantities from a ring of immediate cell neighbors. The resulting stencils involve cells both in the upstream and downstream direction. This means that the resulting discretization graph does not reflect the same unidirectional properties as the underlying flow equations for co-current flow.

Alternatively, one can use discontinuous Galerkin methods (dG), which were first extended to general systems of hyperbolic conservation laws by \citet{Cockburn1989, Cockburn1991}. When combined with a single-point upstream mobility scheme for flux evaluation, these methods preserve the causality of the continuous flow equations. That is, the corresponding stencil is restricted to cells in the upstream direction in regions of co-current flow. After application of flux-based reordering, the discrete nonlinear flow equations are permuted to block-triangular form, with small blocks representing individual cells from regions of co-current flow and larger blocks representing regions of counter-current flow. This gives a natural localization: You start at the inflow locations and solve the nonlinear equations block-by-block toward the location of outflow. For blocks consisting of multiple cells, you can either solve for all mutually dependent unknowns simultaneously, or you can use an effective Gauss--Seidel solver that decomposes the problem to an iteration over a sequence of single-cell problems \citep{Lie2014}. The advantage of the Gauss--Seidel approach is that it relies entirely on single-cell nonlinear solvers, which are simpler to optimize.

In this work, we combine intercell flux reordering methods with higher-order dG methods. This idea has previously been studied by \citet{Natvig07} and \citet{Eikemo09} for passive advection problems and by \citet{Natvig2008} for incompressible problems in two-phase, three-phase, and two-phase--three-component flow without gravity and capillary forces. Later, the first-order variant of the method was extended to include gravity \citep{Lie2012} and to polymer flooding \citep{Lie2014} with compressibility and gravity effects. Herein, we study the widely used family of black-oil equations and present a method that can handle all relevant flow effects including hysteresis and capillary forces as seen in the simulation model of the Norne oil field \cite{norne-data}. We also discuss how to formulate high-order dG discretization for general polyhedral grids, introduce a simple order-reduction method to prevent creating spurious oscillations, and present a new method based on blocks of cells, which is more parallel and cache efficient than solving for a single cell at a time.

\section{Governing equations}

The black-oil model describes conservation of three pseudo-components (water, oil, and gas), which at reservoir conditions can distribute in three phases (aqueous, oleic, and gaseous):
\begin{equation}
\label{eq:black-oil}
    \begin{aligned}
        & \partial_t\left(\phi b_w S_w\right) + \nabla \cdot \left(b_w \vec v_w \right) - b_wq_w = 0, \\
        & \partial_t\left(\phi \left[b_oS_o +  b_g r_{v}S_g\right]\right) & \\
        & \,\, + \nabla \cdot \left(b_o \vec v_o  + b_gr_v \vec v_g\right) - (b_oq_o + b_gr_vq_g) = 0, \\
        & \partial_t\left(\phi \left[b_gS_g +  b_o r_sS_o\right]\right) & \\
        & \,\, + \nabla \cdot \left(b_g \vec v_g  + b_or_s \vec v_o\right) - (b_gq_g + b_or_sq_o) = 0.
    \end{aligned}
\end{equation}
Here, $S_\alpha$ denotes saturation, $\vec v_\alpha$ the macroscopic Darcy velocity, and $q_\alpha$ sources and sinks of phase $\alpha$. Shrinkage factors $b_\alpha$ model pressure-dependent densities $\rho_\alpha$, and the gas-oil and oil-gas ratios $r_s$ and $r_v$ model the volume of gas dissolved in oil and oil vaporized in gas, respectively, both at standard conditions. The phase velocity $\vec v_\alpha$ is given by the multiphase extension of Darcy's law:
\begin{equation*}
    \label{eq:darcy-velocity}
    \vec v_\alpha = - \lambda_\alpha \tens K \left(\nabla p_\alpha - \rho_\alpha g \nabla z\right),
\end{equation*}
where $\lambda_\alpha = k_{r\alpha}/\mu_\alpha$ is the mobility of phase $\alpha$; relative permeability $k_{r\alpha}$ models the reduced mobility of one phase in the presence of another, whereas $\mu_\alpha$ as usual denotes the viscosity of phase $\alpha$. The model is closed by assuming that the fluid phases fill up the entire pore space, $S_w + S_o + S_g = 1$, and that the phase pressures are related through capillary pressures:
\begin{equation*}
    \label{eq:capillary-pressure}
    p_o = p_w + P_{cow}(S_w, S_o), \quad p_g = p_o + P_{cgo}(S_o, S_g).
\end{equation*}

\section{Sequential splitting}

We use backward Euler to discretize in time. On semi-discrete form, the equation describing conservation of mass for water \eqref{eq:black-oil} at time $n+1$ then reads
\begin{equation}
\label{eq:water-semidiscrete}
    \begin{aligned}
        \res_w = & \frac{1}{\Delta t} \left( \left[\phi b_w S_w\right]^{n+1} - \left[\phi b_w S_w\right]^{n}\right) \\
        & \qquad + \nabla \cdot \left(b_w \vec v_w \right)^{n+1} - (b_wq_w)^{n+1} = 0,
    \end{aligned}
\end{equation}
where we have introduced the time-step length $\Delta t = t^{n+1} - t^n$, and superscripts $n$ and $n+1$ refer to the time step. Semi-discrete forms for the oil ($\res_o$) and gas ($\res_g$) equations are analogous. We multiply each equation by the following factors
\begin{align*}
  \omega_w & = \frac{1}{b_w^{n+1}}, \\
  \omega_o & = \frac{1}{1-r_s^{n+1}r_v^{n+1}}\left(\frac{1}{b_o^{n+1}} - \frac{r_{s}^{n+1}}{b_g^{n+1}}\right), \\
  \omega_g & = \frac{1}{1-r_s^{n+1}r_v^{n+1}}\left(\frac{1}{b_g^{n+1}} - \frac{r_{v}^{n+1}}{b_o^{n+1}}\right),
\end{align*}
and sum to eliminate all terms involving saturations at time step $n+1$. This gives us a pressure equation 
\begin{equation}  
  \label{eq:pressure-semidiscrete}
  \res_p = \omega_w \res_w + \omega_o \res_o + \omega_g \res_g.
\end{equation}
To obtain a fully discrete formulation, we introduce a grid covering our computational domain with $\ncell$ non-overlapping polyhedral cells $\{\cell_i\}_{i = 1}^{\ncell}$. We denote the common interface of cell $i$ and $j$ by $\face_{ij}$ and the volumetric flux from cell $i$ to cell $j$ by
\begin{equation*}
  v_{\alpha,ij} = \int_{\face_{ij}} \vec v_\alpha \cdot \n_{ij} \d \sigma.
\end{equation*}
Here, $\vec n_{ij}$ is the unit normal from cell $i$ to cell $j$. Notice that conservation of mass requires that $v_{\alpha, ij}=-v_{\alpha, ji}$. In the following, we assume that we have a pressure solver that can solve the fully discrete version of \eqref{eq:pressure-semidiscrete} to yield a total velocity field $\vec v = \vec v_w + \vec v_o + \vec v_g$, given as a set of constant fluxes $v_{ij} = v_{w, ij} + v_{o, ij} + v_{g, ij}$ over each cell interface. In the numerical experiments, unless otherwise noted, we use the standard black-oil pressure solver from the open-source MATLAB Reservoir Simulation Toolbox \citep{Lie2016MRST-book}.

In a sequential solution procedure, we solve \eqref{eq:pressure-semidiscrete} with fixed saturations to obtain pressures and the total Darcy velocity $\vec{v}$. Given the total flux, we compute new phase fluxes using a fractional flow formulation. These are defined from the following formula in the semi-continuous case (neglecting capillary effects)
\begin{align}
    \label{eq:phase-velocity}
    \vec v_{\alpha} = f_{\alpha}\left( \vec v + \tens K \sum_{\beta \neq \alpha} \lambda_{\beta} \left(\rho_\alpha - \rho_\beta\right) g \nabla z\right),
\end{align}
where we have introduced the fractional flow function
\begin{equation*}
    f_{\alpha}(S_w, S_o, S_g) = \frac{\lambda_{\alpha}(S_\alpha)}{\lambda_{w}(S_w) + \lambda_{o}(S_o) + \lambda_{g}(S_g)}.
\end{equation*}
Capillary pressure differences are accounted for with a term similar to the gravity term in \eqref{eq:phase-velocity}. Using this, we solve the transport equations \eqref{eq:black-oil} with fixed pressure to obtain new saturations (or solution ratios). This introduces a splitting error proportional to the time step $\Delta t$ that can be significant for cases with strong coupling between pressure and saturation. One can also make the solution converge towards the fully implicit solution by adding outer iterations, as described by e.g, \citet{Jenny2006}.

\section{Discontinuous Galerkin discretization}

For simplicity, we only describe the weak formulation for the water equation without capillary and gravity effects, even though gravity effects are included in our solver. Capillary effects are also included, but currently only supported in our first-order transport solver. Moreover, we drop the phase subscript $w$ and the time superscript $n+1$. Equation \eqref{eq:water-semidiscrete} now reads
\begin{equation}
    \label{eq:water-residual}
    \res_w = \tfrac{1}{\Delta t} \left( \left[\phi b S\right] - \left[\phi b S\right]^{n}\right) + \nabla \cdot \left(b f\vec v \right) - (bq) = 0,
\end{equation}
where we recall that the fractional flow function $f$ depends on all phase saturations. We multiply by a test function $\psi$ in a function space $V$ of arbitrarily smooth functions and integrate (by parts) over cell $\cell_i$:
\begin{equation*}
    \label{eq:water-weak}
    \begin{aligned}
        \frac{1}{\Delta t} \int_{\cell_i} \Big( (\phi b S)- (\phi b S)^{n}\Big) \psi \d V - &\int_{\cell_i} b f \vec v \cdot \nabla\psi \d V \\
        + \int_{\partial \cell_i}  \left[b f \psi \right] \vec v \cdot \vec n \d \sigma - &\int_{\cell_i} (bq)\psi \d V= 0.
    \end{aligned}
\end{equation*}
Here, the square brackets in the surface integral signify that the integrand is possibly discontinuous across the cell boundary. Note that we can express the surface integral as
\begin{equation*}
  \int_{\partial \cell_i} \left[b f \psi \right] \vec v \cdot \vec n \d \sigma = \sum_{j\in \neigh (i)} \int_{\face_{ij}} \left[b f \psi \right] \vec v \cdot \vec n_{ij} \d \sigma,
\end{equation*}
where $\neigh(i)$ denotes the set of all cells sharing a common interface with cell $i$. We define the upwind operator 
\begin{equation*}
    \begin{aligned}
        \upw[y, \vec v_\alpha \cdot \vec n_{ij}; x] & =
        \begin{cases}
            y_{ij} \vec v_\alpha \cdot \vec n_{ij}, & \text{if} \quad \vec v_\alpha \cdot \vec n_{ij} > 0, \\
            y_{ji} \vec v_\alpha \cdot \vec n_{ij}, & \text{if} \quad \vec v_\alpha \cdot \vec n_{ij} \leq 0,
        \end{cases}\\
        \text{where } y_{ij} & = \lim_{x'\to x} y(x')
        \text{ for } x'\in\cell_i.
    \end{aligned}
\end{equation*}
Subscripts $ij$ and $ji$ denote the limiting interface value of $y(x')$ for $x \in \face_{ij}$ as $x'$ approaches $x$ from inside cell $i$ or cell $j$, respectively. If the coordinate is not important, we simply write $\upw[y, \vec v_\alpha \cdot \vec n_{ij}]$. Using this, we write
\begin{equation*}
    \int_{\partial \cell_i}\left[b f \psi \right] \vec v \cdot \vec n \d \sigma = \sum_{j\in N(i)}  \int_{\face_{ij}} \upw\left[b, f\vec v \cdot \vec n_{ij}\right] \psi \d \sigma.
\end{equation*}
Since the fractional flow function is always non-negative, $\vec v$ and $\vec v_\alpha = f_\alpha \vec v$ will point in the same direction in the absence of  gravity and capillary effects, and the upwind definition is explicit. However,
referring back to \eqref{eq:phase-velocity}, we see that the upwind definition is generally implicit in the presence of gravity/capillary effects, since quantities in the parenthesis depend on the mobilities $\lambda_\alpha$. An equivalent explicit upwind definition has been derived by \citet{brenier1991upstream} for finite-volume methods. Herein, we employ this definition \textit{at each cubature point} involved in evaluating the surface integrals.

With this notation, we define the following weak form of the residual equation \eqref{eq:water-residual} in cell $i$
\begin{equation}
    \label{eq:water-bilinear}
    \begin{aligned}
    &\res_{w,i}(S_w, S_o, S_g, \psi) \\
    & = \tfrac{1}{\Delta t} \acc_{w,i}(S_w, \psi) + \flux_{w,i}(S_w, S_o, S_g, \psi) - \src_{w,i}(\psi),
    \end{aligned}
\end{equation}
where the accumulation, flux and source/sink functionals are defined as
\begin{align*}
    \acc_{w,i}(S_w, \psi) = &\int_{\cell_i} \Big( (\phi b_w S_w)- (\phi b_w S_w)^{n}\Big)\psi \d V,\\
        \flux_{w,i}(S_w,S_g,S&_o,\psi) \\
         = \sum_{j\in N(i)}& \int_{\face_{ij}} \upw\left[b_w, f_w(S_w,S_o,S_g)\vec v \cdot \vec n_{ij}\right] \psi \d \sigma \\
         -&\int_{\cell_i} b_w f_w(S_w,S_o,S_g) \vec v \cdot \nabla \psi \d V,\\
    \src_{w,i}(\psi) = &\int_{\cell_i} b_wq_w \psi \d V.
\end{align*}
Weak form residuals for the oil ($\res_{o,i}$) and gas ($\res_{g,i}$) can be derived in a similar fashion. To obtain a discrete weak formulation, we replace the function space $V$ with a finite-dimensional subspace $V_h$ consisting of functions that are smooth on each cell $\cell_i$, but possibly discontinuous across cell interfaces. We replace the saturations and the test function $\psi$ with approximations $S_{\alpha,h}\in V_h$ and $\psi_h\in V_h$, and arrive at the following discrete weak formulation:
\begin{equation*}
    \label{eq:water-weak-disc}
    \begin{aligned}
        \text{Find } S_{w,h},&S_{o,h},S_{g,h} \in V_h \text{ such that }\\ 
        \res_{\alpha,i}(&S_{w,h},S_{o,h},S_{g,h}, \psi_h) = 0 \text{ for all } \psi_h \in V_h.
    \end{aligned}
\end{equation*}
If we introduce the basis $\{\psi_k\}_{k = 1}^{\ndof}$ for $V_h$, we may express the saturations as $S_{\alpha,h} = \sum_{k=1}^{\ndof} S_{\alpha,k} \psi_k$, where $S_{\alpha,k} \in \real$ is referred to as the $k$th degree of freedom of $S_\alpha$. The discrete weak formulation then takes the form
\begin{equation*}
    \label{eq:water-weak-mat}
    \begin{aligned}
        & \text{Find } (S_{\alpha,1}, \dots, S_{\alpha,\ndof}) \in \real^{\ndof}, \, \alpha = w,o,g, \text{ such that } \\
        & \,\res_{\alpha,i}\left(\sum_{k=1}^{\ndof} S_{w,k}\psi_k ,\sum_{k=1}^{\ndof} S_{o,k}\psi_k,\sum_{k=1}^{\ndof} S_{g,k}\psi_k, \psi_\ell\right) = 0 \\
        & \,\text{for all } \ell = 1, \dots, \ndof.
    \end{aligned}
\end{equation*}
Note that we typically use the closure relation $S_w + S_o + S_g = 1$ to eliminate one variable and one equation from the above.

There are several important factors that need careful consideration to implement the above discretization. In the following, we briefly discuss the most important of these. A more detailed analysis is subject of further research.

\subsection{Basis functions}

\begin{figure}[t]
    \centering
    \includegraphics[width = 0.48\textwidth]{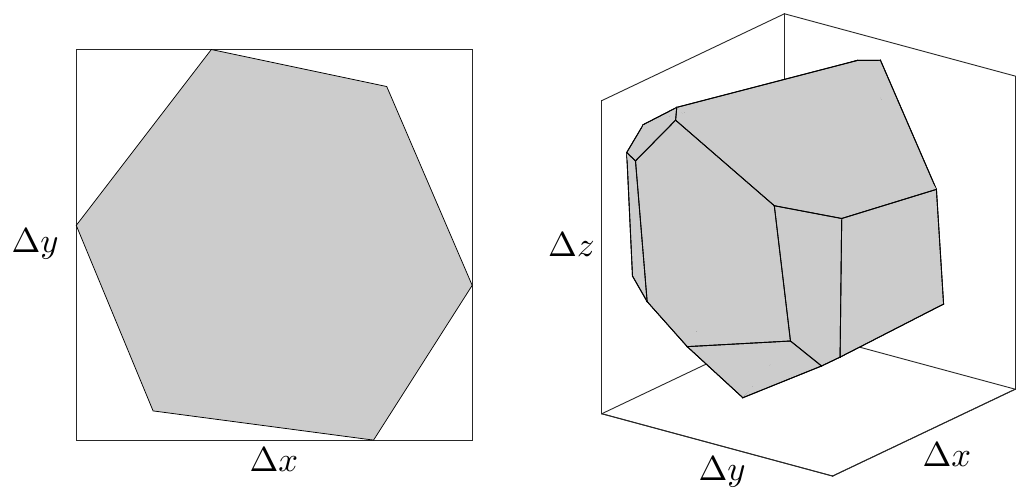}
    \caption{Dimensions of bounding boxes used to define the dG basis functions for a polygon in 2D and a polyhedron in 3D.}
    \label{fig:basis-functions}
\end{figure}

In finite-element methods, it is common to use \textit{orthogonal} basis functions $\{\psi_k\}_{k = 1}^{\ndof}$, i.e., $\int_{\Omega_i} \psi_k \psi_\ell = \delta_{k,\ell}$. The main reason for this choice is that it typically results in a significantly less dense discretization matrix. Defining orthogonal basis functions for general polyhedral grids is not a trivial problem. For the types of problems we are interested in, however, the weak form residual \eqref{eq:water-bilinear} will generally depend nonlinearly on the saturations (and thus the basis functions) due to the fractional flow function $f_\alpha$ in the flux terms $\flux_{\alpha,i}$. Effectively, orthogonal basis functions will generally not lead to less dense discretization matrices. Herein, we will therefore simply use tensor products of Legendre polynomials as our basis functions. These form a convenient basis for the space of polynomials and are orthogonal for cuboid grid cells. The first two Legendre basis functions are $\ell_0(x) = 1$ and $\ell_1(x) = x$, and higher-order basis functions are defined in our dG implementation by exploiting Bonnet's recursion formula,
\begin{equation*}
    (k+1)\ell_{k+1}(x) - (2k+1)x\ell_k(x) + k\ell_{k-1}(x) = 0.
\end{equation*}
Basis function $j$ of polynomial order $r+s+t$ for cell $i$ is then defined as
\begin{equation*}
    \label{eq:legendre-basis}
    \psi_{j}(\vec x) =
    \begin{cases}
        \ell_r\left(\frac{x - x_i}{\Delta x_i/2}\right)\ell_s\left(\frac{y - y_i}{\Delta y_i/2}\right)\ell_t\left(\frac{z - z_i}{\Delta z_i/2}\right), & \vec x \in \cell_i, \\
        0, & \vec x \not\in \cell_i,
    \end{cases}
\end{equation*}
where
\begin{equation*}
    \vec x = (x, y, z), \quad \vec x_i = (x_i, y_i, z_i), \quad \Delta \vec x = (\Delta x_i, \Delta y_i, \Delta z_i)
\end{equation*}
denote the spatial coordinate, cell centroid and the dimensions of the smallest cuboid aligned with the coordinate axes that completely contains cell $i$, as shown in Figure~\ref{fig:basis-functions}. For a dG method of degree $k$, the basis consists of all tensor products on this form such that $  0 \leq r+s+t \leq k$, with $r,s,t\geq 0$. In $d$ space dimensions, this gives a total number of basis functions
\begin{equation*}
  \ndof = \binom{k+d}{d} = \frac{(k+d)!}{d!k!}.
\end{equation*}
We use the notation dG($k$) to refer to a dG method of degree $k$. The error of a method of degree $k$ will for smooth solutions decay with order $k+1$, so that a dG($k$) method is of formal order $k+1$.

\subsection{Cubature rules}

The integrals in \eqref{eq:water-bilinear} must be numerically evaluated. Efficient cubature rules for general polyhedral grid cells are crucial to obtain an efficient implementation. A straightforward approach to numerically evaluate integrals over polyhedral cells is to partition the cell into non-overlapping simplices and apply a known cubature rule for these. This uses significantly more cubature points than strictly needed for the cubature to be correct and is thus highly inefficient from a computational point of view. Instead, we apply an approach based on moment-fitting, see e.g., \cite{Muller2013}. This approach defines the quadrature rule for a cell $\Omega$ (omitting subscript $i$) as the solution to a small linear system
\begin{align*}
    & \begin{bmatrix}
        \psi_{1}(\vec x_{1}) & \dots & \psi_{1}(\vec x_{\ndof}) \\
        \vdots  & \ddots & \vdots \\
        \psi_{\ndof}(\vec x_{1}) & \dots & \psi_{\ndof}(\vec x_{\ndof}) \\
    \end{bmatrix}
    \begin{bmatrix}
        w_{1} \\ \vdots \\ w_{\ndof}
    \end{bmatrix}
    \\ & \hspace{6em} =
    \frac{1}{|\Omega|}\begin{bmatrix}
        \int_{\cell} \psi_1 \d V\\
        \vdots \\
        \int_{\cell} \psi_{\ndof} \d V
    \end{bmatrix}, \quad
    \text{or} \quad
    \Psi \vec w = \vec b.
\end{align*} 
The integrals (or moments) on the right-hand side can be evaluated using a sub-optimal cubature rule. If we pick the cubature points $\vec x_0, \dots \vec x_{\ndof}$ so that the matrix $\Psi$ is invertible, we can solve for the quadrature weights $w_1, \dots w_{\ndof}$ to obtain a cubature rule with $\ndof$ points for each cell. Note that it is possible to construct a quadrature rule of precision $k$ with less than $\ndof$ points by eliminating points with marginal significance. This computation can be done once for each grid cell in a preprocessing step.

In 3D, the surface integrals in the second term of $\flux$ in \eqref{eq:water-bilinear} are evaluated in the same way as areal integrals in 2D, using the moment-fitting approach described above. In 2D, the surface integrals are simple line integrals, for which a standard Gauss cubature rule is employed.

\subsection{Velocity interpolation}

As stated above, we assume that the solution to the pressure equation yields a total velocity as a set of fluxes that are constant on each interface. We see from \eqref{eq:water-bilinear} that for linear and higher-order basis functions, we need information about the velocity $\vec v$ also inside the cell, meaning that we must interpolate from the interface fluxes. Herein, we apply a simple scheme inspired by the mimetic finite difference method \citep{MRST12:comg}. This gives a constant velocity inside the cell, consistent with the volumetric interface fluxes $v_{ij}$:
\begin{equation*}
  \vec v_i = \frac{1}{|\cell_i|}\sum_{j \in N(i)} v_{ij}(\vec x_{ij} - \vec x_i).
\end{equation*}
Here, $\vec x_{i}$ and $\vec x_{ij}$ refer to the centroid of cell $i$ and interface $ij$, respectively. The velocity on an interface is assumed to be constant. More sophisticated interpolation schemes like the extended corner-velocity interpolation scheme (ECVI) \citep{Klausen2012} have earlier been investigated for dG methods for flow diagnostics \citep{Rasmussen2014}, but will for simplicity not be applied herein.

\subsection{Order reduction}
Higher-order methods usually require a slope limiter, or some other means to prevent the creation of over- and under-shoots or other types of nonphysical oscillations near spatial discontinuities in the solution. \citet{Natvig07} previously showed that order reduction followed by local and dynamic grid refinement could be used as an alternative strategy to avoid oscillations and give high-resolution in a robust and cost-efficient manner. Herein, we use a simplified version: Whenever the jump across an intercell interface exceeds a prescribed tolerance, we reduce the local approximation space to dG(0). The same is done when the solution exceeds its domain of definition by a small factor $\varepsilon$.   

\section{Localized nonlinear solvers}

This section discusses how we can utilize certain unidirectional flow properties to formulate efficient nonlinear solution strategies that localize the Newton procedure to avoid spending unnecessary iterations in regions of the reservoir where many iterations are not needed.

\subsection{Reordering based on intercell fluxes}

An important physical property of the transport \eqref{eq:black-oil} is that the flow is unidirectional along streamlines if no gravity or capillary effects are present. Each dG basis function is restricted to a single cell only, and hence the only intercell interactions in the weak form $\res_{w,i}$ in \eqref{eq:water-bilinear} are due to the upwind operator. We split the neighboring cells $\neigh(i)$ of $i$ into upstream cells $\upstr(i)$ and downstream cells $\downstr(i)$ based on the total flux $\vec{v}\cdot\vec{n}_{ij}$:
\begin{align*}
    \sum_{j\in \neigh(i)} & \int_{\face_{ij}} \upw\left[b, f\vec v \cdot \vec n_{ij}\right] \psi \d \sigma &\\
    & = \sum_{j\in \upstr(i)} \int_{\face_{ij}} (b f)_j\, \vec v \cdot \vec n_{ij} \psi \d \sigma & (\text{Upstream})\\
    & + \sum_{j\in \downstr(i)} \int_{\face_{ij}} (b f)_i\, \vec v \cdot \vec n_{ij} \psi \d \sigma, & (\text{Downstream})
\end{align*}
where $(\cdot)_i$ indicates that quantities inside the parentheses should be evaluated from cell $i$. Recall that all other unknowns in the weak-form residual $\res_{w,i}$ have support limited to cell $i$, and hence the inherent unidirectional property of the continuous equations \eqref{eq:black-oil} is preserved by the discretization. In practice, this means that if we know the solution in all upstream cells $\upstr(i)$, the only unknowns in $\res_{w,i}$ are the ones associated with cell $i$. That is, if we now perform a topological sort of the directed, acyclic graph (DAG) induced by the intercell fluxes and use this to reorder the cells, we can solve the transport equations cell-by-cell by traversing the sorted graph. The algebraic interpretation of this is the following: If we linearize the nonlinear transport equations for all cells simultaneously, and permute the system according to the topological order, we obtain a lower-triangular matrix. Note, however, that we never assemble the discretization matrix for the full system in the reordering solution procedure. Instead we solve the nonlinear transport equations $\res_{\alpha,i} = 0$ cell-by-cell. This way, we avoid expensive linearizations of large systems of nonlinear equations.

\begin{figure*}[t!]
    \centering
    \includegraphics[width = \textwidth]{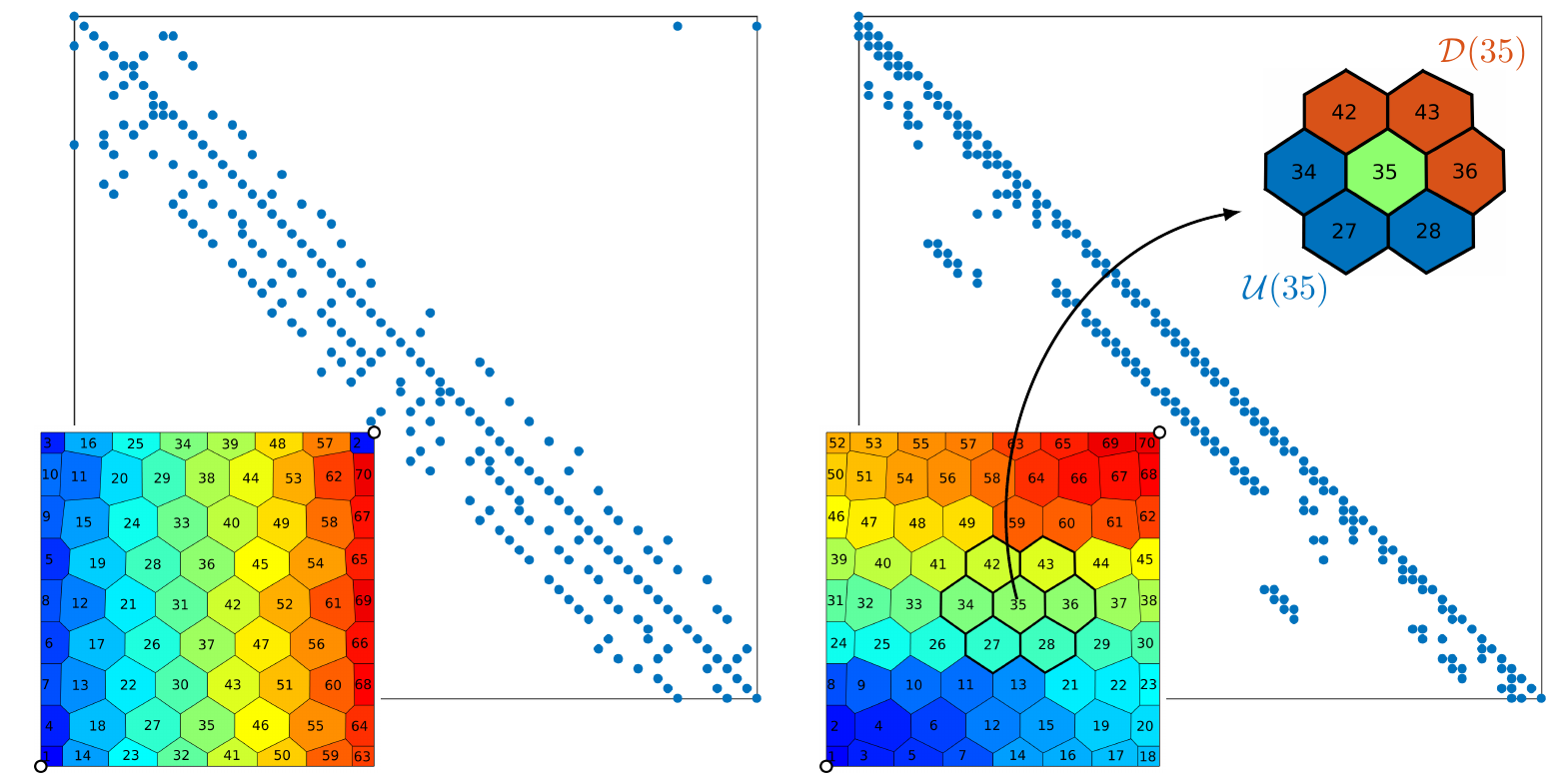}
    \caption{Sparsity pattern for a dG(0) method formulated on a Voronoi grid with original (left) and reordered (right) numbering of cells. The inset shows grid cell 35 and its upstream $\upstr(35)$ and downstream $\downstr(35)$ neighbors.}
    \label{fig:reordering-dG0}
\end{figure*}

Figure~\ref{fig:reordering-dG0} illustrates the reordering principle applied to a dG(0) scheme for a quarter-five spot test case formulated on a 2D Voronoi grid. In the original grid, the cells are ordered almost, but not entirely, by the $x$-coordinate of their centroids, giving an irregular sparsity pattern with entries both above and below the diagonal. After reordering, the cells are ordered such that if $j\in \neigh(i)$, then $j\in \upstr(i)$ for $j<i$ and $j\in \downstr(i)$ for $j>i$, where we recall that $\neigh(i)$ denotes neighbors and $\upstr(i)$ and $\downstr(i)$ are cells that lie upstream/downstream of cell $i$.

\begin{figure}[b]
    \centering
    \includegraphics[width=\linewidth]{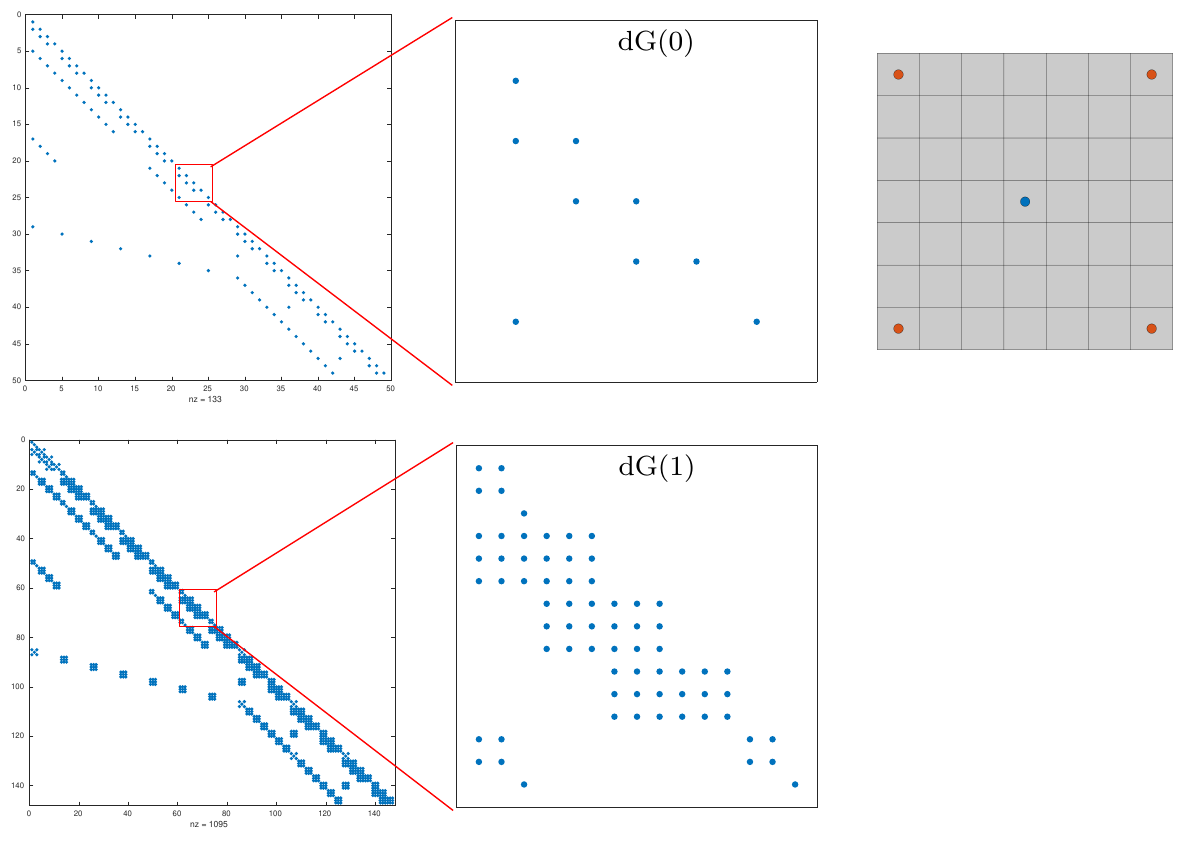}
    \caption{Sparsity pattern after reordering for dG(0)/dG(1) for a five-spot problem on a rectangular grid.}
    \label{fig:reordering-dG1}
\end{figure}

Figure~\ref{fig:reordering-dG1} similarly illustrates the difference between a dG(0) and a dG(1) discretization on a Cartesian grid. The two matrices have exactly the same flux graph, and thus obtain the same cell ordering. However, whereas dG(0) gives a nonlinear scalar problem in each cell, dG(1) gives a $3\times 3$ nonlinear system.

\begin{figure*}[t!]
    \centering
    \includegraphics[width=\linewidth]{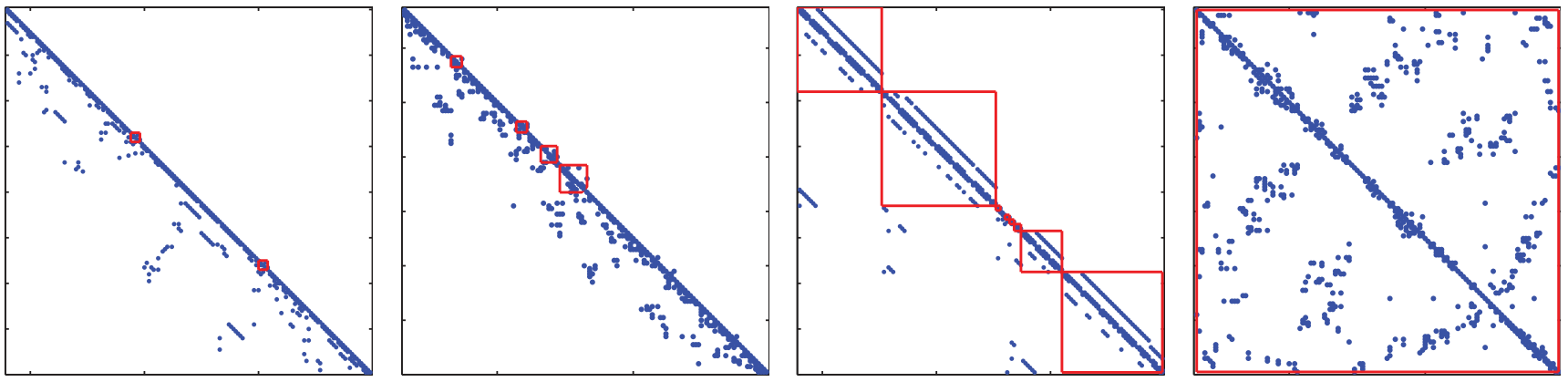}
    \caption{Examples of different degrees of coupling, where the cycles can be seen as blocks on the diagonal of the dG(0) discretization matrix. In the first sparsity pattern, the system is almost completely reorderable with only two small cycles; in the second, a small fraction of the cells belongs to four small cycles; the third patterns shows an example where almost all cells belongs to eight larger cycles; finally the last pattern shows a case where all cells belong to one huge cycle. All cells in the same cycle must be solved simultaneously. Figure from \cite{Lie2014}.}
    \label{fig:cycles}
\end{figure*}

When gravity and capillary effects are included, the flow will no longer be unidirectional along streamlines. As a result, the intercell flux graph may contain cycles, so that the reordering method cannot be applied directly. However, we may construct a DAG from the intercell flux graph by grouping cells that are part of the same cycle into a single supernode. By topologically sorting this DAG, we end up with a block-triangular system in which each block consists of degrees of freedom that are mutually coupled and must be solved for simultaneously. Note that as explained above, we evaluate the upstream direction for each quadrature point at the interface. For a higher-order method, the upstream direction may thus vary along a single interface when gravity/capillary effects are included, particularly when we use high polynomial order and/or cells with large vertical extent. This means that each quadrature point on interfaces with changing upstream direction must be treated as an edge in the original intercell flux graph, and such interfaces will automatically introduce a local cycle and require the neighboring cells to be grouped into a supernode. In general, the reordering procedure will result in a DAG with $\ncomp$ nodes, where $\ncomp$ equals the number of connected components of the original intercell flux graph. Topologically sorting this graph gives us an ordering
\begin{equation*}
    \order = \left(\order_1, \dots, \order_{\ncomp}\right), \quad \order_k = (i_1, \dots, i_{n_k}),
\end{equation*}
where each $\order_k$ is a set of cell indices, and $n_k$ is the number of cells in connected component $k$. The number of connected components $\ncomp$ may vary from one in the worst case, where all cells are connected in a single cycle, to the total number of cells $\ncell$ in the best case, where the problem can be solved cell by cell without the need of any supernodes. Given two connected components $\order_k$ and $\order_\ell$, we will then have that if $i>j$ for a cell $i \in \order_k$ and a cell $j \in \order_\ell$ then $i>j$ for all cells $i \in \order_k$, $j \in \order_\ell$. Figure~\ref{fig:cycles} illustrates different scenarios of coupling, from an almost completely reordered system, to a system where all cells belong to a single cycle. 

\subsection{Block-wise processing}

For test cases with only viscous forces, the reordered system can be solved cell-by-cell in a single pass. This is optimal in the sense that it minimizes the number of cell-wise nonlinear solves over the simulation. However, it is not necessarily optimal in terms of runtime; most modern workstations contain several CPU cores, as well as any number of SIMD-type registers. In this setting, feeding sufficient computations to the CPU is often a large bottleneck and a barrier to rapid execution time. For this reason, our solver has an optional block-wise algorithm that solves a prescribed number of cells simultaneously. If the largest ordering index of the previous solve was $i$, this algorithm then solves for indices $i+1$ to $i+n_b$, where $n_b$ is the block size. This gives a larger system, and individual cells can follow a sub-optimal solution path when solved by a multi-cell Newton--Raphson solver. On the other hand, the startup cost is incurred only once per block of cells and not for every cell, and the block solve utilizes cache better and leaves less computational resources idle. Hence, simultaneous solve of $n_b$ cells will typically be significantly faster than $n_b$ single-cell solvers.  

\begin{figure*}[t!]
    \centering
    \includegraphics[width=1\linewidth]{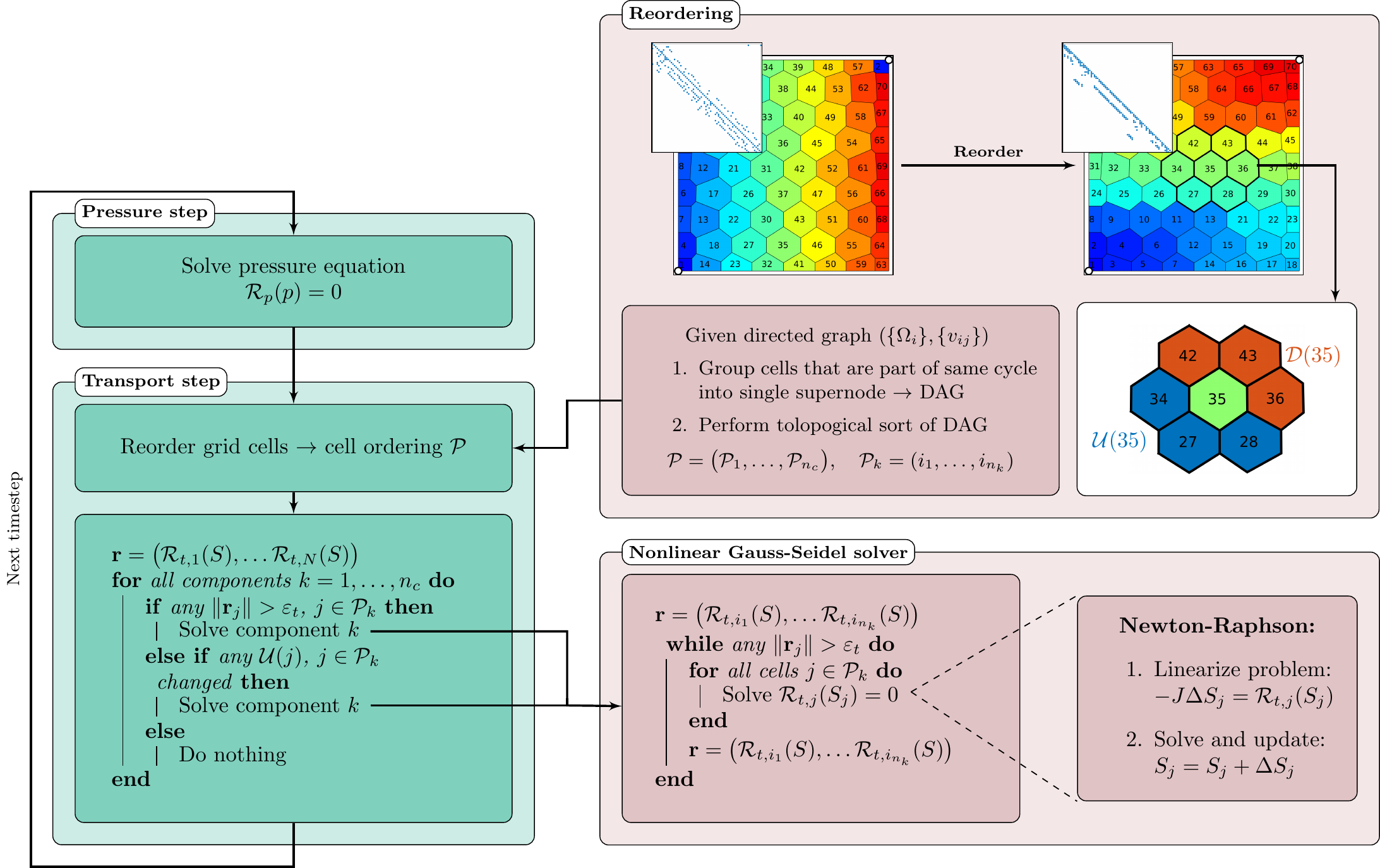}
    \caption{Schematic overview of the solution procedure. The pressure equation is solved with a suitable solver, resulting in intercell fluxes $v_{ij}$, which we may interpret as a directed graph with grid cells as nodes, and intercell fluxes as edge weights. By grouping all cells that are part of the same cycle into a (super)node, we obtain a directed, acyclic graph (DAG). A topological sort of this graph gives us a node (or component) ordering $\order$. By traversing the nodes in this order, we can solve the transport problem $\res_t(S) = 0$ node by node, since the cells in each node in the DAG only depend on each other and cells from upstream (and already resolved) nodes. }
    \label{fig:solution-procedure}
\end{figure*}

\subsection{Nonlinear Gauss--Seidel solvers}
As explained above, gravity will result in a block-triangular system in which each block consists of degrees of freedom that are mutually coupled and must be solved for simultaneously. The size of each block depends on whether the block represents degrees of freedom that are part of a co-current or counter-current flow region. Blocks in co-current regions represent the discrete system posed inside a single cell and are treated as described above. Blocks in counter-current regions consist of multiple cells that are coupled in cycles. In principle, each of these cycles of cells should be solved by a localized Newton iteration. Another possibility is to first solve the water equations cell-by-cell in the order of descending water phase potential, and then solve the oil equations in the order of descending oil phase potential \cite{Kwok2007}. This method requires that the oil component equation only depends on water and oil saturation, which is the case for black-oil models. Herein, we use a nonlinear Gauss--Seidel strategy instead that simply loops through all the cells of a cycle and solves them one at a time, using the values in all neighboring cells as if they were all known upstream values. This will typically not give the correct solution after a single iteration through the loop, but a large number of numerical experiments indicate that the process converges after only a few repetitions. Capillary forces will in the worst case couple all cells in one big cycle, as seen in the rightmost sparsity pattern in Figure~\ref{fig:cycles}. However, the coupling is oftentimes so weak that a Gauss--Seidel type iteration is still very effective.

\section{Solution procedure}
The entire procedure to evolve the solution one time step can now be summarized as follows: (\textit{i}) Solve the pressure equation \eqref{eq:pressure-semidiscrete} to obtain intercell fluxes $v_{ij}$ and oil phase pressure $p_o$; (\textit{ii}) construct a DAG from the intercell flux graph by grouping cells that are part of the same cycle into a single supernode; (\textit{iii}) starting from the first node in the sorted graph, solve the transport equations \eqref{eq:water-bilinear} cell-by-cell or cycle-by-cycle in topological order. The procedure is summarized in Figure~\ref{fig:solution-procedure}. Notice the important property that if the residual of a cell or cycle is converged at the beginning of the time step, and none of the upstream neighbors were updated in this step, we do not need to call the nonlinear Gauss--Seidel solver and can proceed to the next component. Effectively, we thus avoid solving for a large number of zeros, and can instead focus our computational resources on parts of the domain where updates are nonzero.

\begin{figure*}[t]
    \centering
    \includegraphics[width = \textwidth]{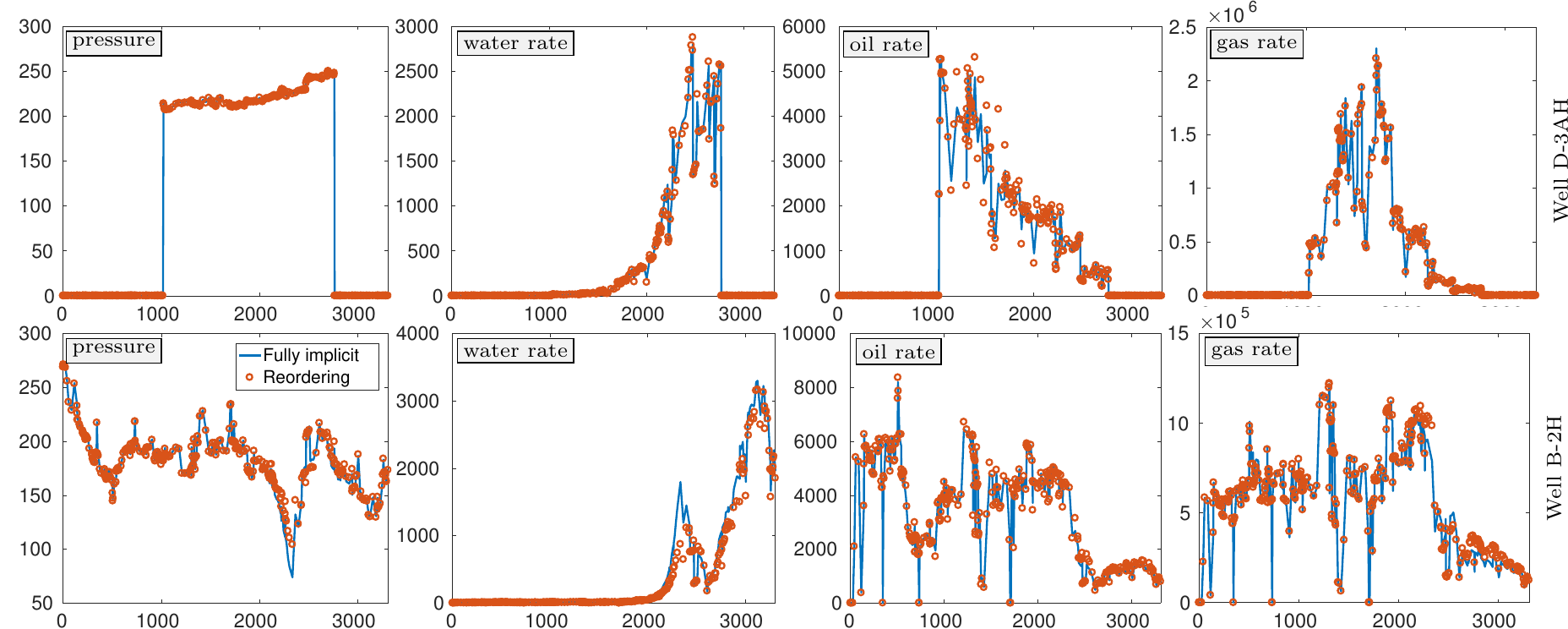}
    \caption{Simulation responses for two of the many wells in the Norne field model. We plot the bottom-hole-pressure and production rates for water, oil and gas at surface conditions.}
    \label{fig:norne-well}
\end{figure*}

\section{Numerical examples}
We have implemented our reordering approach in two different open-source simulator platforms. The MATLAB Reservoir Simulation Toolbox (MRST) \citep{Lie2016MRST-book}  is a general toolbox aimed at rapid development of proof-of-concept simulators and workflow tools. We have used MRST to develop a combined dG--reordering method that allows for cell-based refinement. The solver uses automatic differentiation to linearize the localized systems and is as such applicable to a wide variety of multiphase flow models. At the time of writing, the implementation in MRST has only been tested for two-phase models with varying degrees of compressibility. 

OPM Flow \citep{opm} is developed with the goal of providing industry-grade simulations and generally offers single-core computational performance that is comparable with contemporary commercial simulators. We have used this as a platform to implement dG(0) with reordering for a general three-phase black-oil model with compressibility, dissolution and vaporization, gravity and capillary forces, and hysteresis. Our implementation uses a simple, global nonlinear Gauss--Seidel method that traverses all cells repeatedly until the residual of the transport equations are below a prescribed tolerance in all cells.

\subsection{The Norne oil field}
\label{sec:norne}

\begin{figure*}[h!]
    \centering
    \includegraphics[width = 0.95\textwidth]{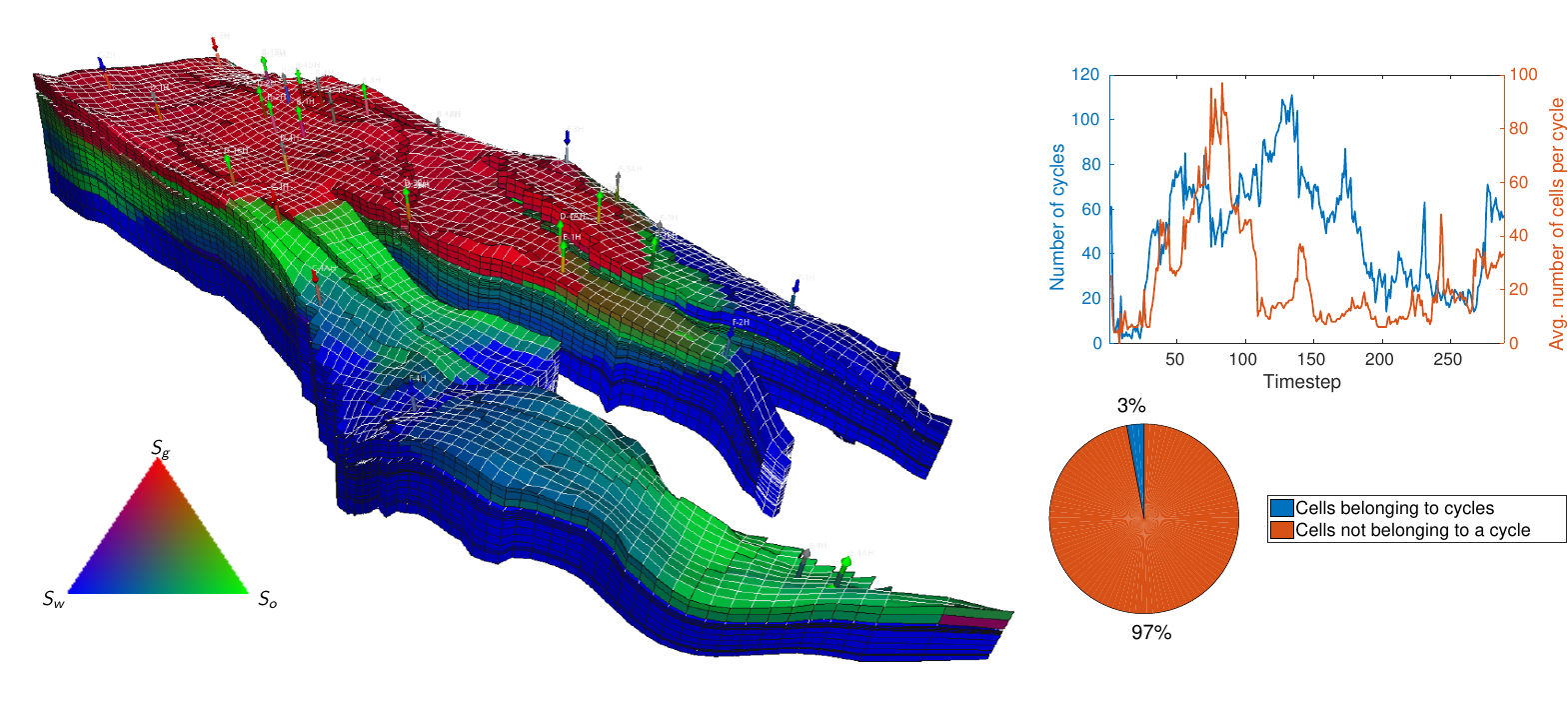}
    \caption{Initial fluid saturations and wells (left); aggregate statistics for the number of cycles (right) for the Norne field model in Example~\ref{sec:norne}.}
    \label{fig:norne-setup}
\end{figure*}

The Norne oil field in the Norwegian Sea is one of the few real assets for which simulation models and other data have been made openly available \cite{norne-data}. Over the 3,312 days of production history, 36 wells are active. We study a water-alternating-gas (WAG) injection scenario, in which the wells are controlled by reservoir-volume rates to match the historical rates to obtain a reasonable pressure development. The purpose of the example is to demonstrate that the sequential splitting method is capable of capturing the correct behavior of an industry-grade simulation run in history-matching mode.

We ran the case with two simulators from OPM: the fully implicit \texttt{Flow} simulator and the experimental \texttt{Flow-reorder} simulator, which uses sequential splitting and reordering of the transport equations. We note that these simulators both support a slightly more complex fluid model than that described in the introduction, incorporating hysteresis effects for relative permeability. Both simulators are available from OPM in source and binary form \citep{opm}.

We highlight the behavior of two different wells to compare the fully implicit to the sequential/reordering solution. Well ``D-3AH'' is a typical producer and shows good match between the fully implicit and reordering solutions for both bottom-hole pressure and rates; see the upper row in Figure~\ref{fig:norne-well}. Well ``B-2H'' is also a producer, which we have chosen to highlight because it shows the worst match between the two solutions. As can be seen in the lower row of Figure~\ref{fig:norne-well}, gas and oil production rates match well, but water-production rate and bottom-hole pressure show significant deviation some time after water breakthrough. There are two main factors that contribute to the difference: First of all, the sequential splitting in \texttt{Flow-reorder} does not use outer iterations. The overall solution algorithm is thus not guaranteed to reduce the overall residual at the end of each time step to the same tolerance as the fully implicit method in \texttt{Flow}. Moreover, convergence is measured somewhat differently in the global Gauss--Seidel and the Newton--Raphson methods.

In our opinion, this simulation is proof of concept for the reordering idea in a production-grade code running a full simulation of a real field in history mode. This is a first step, and there are several ways in which this implementation can be improved. The important observation is that we made no simplifications in the description of geology, fluid behavior, and well and production facilities.

Figure~\ref{fig:norne-setup} shows the fluids in place and reports the occurrence of cycles throughout the simulation. Altogether, the model has approximately 44,000 active cells. The number of cycles and the average number of cells in each cycle vary largely throughout the simulation because of significant changes in well controls and fluid composition. In any time step, the number of cells belonging to cycles only represents a fraction of the total cells. This shows the inherent locality and the speedup potential, since most of the cells can be solved independently.

\subsection{Subset of the SPE10 Model 2}
\label{sec:spe10}

\begin{figure}[t!]
    \centering
    \includegraphics[width = 0.48\textwidth]{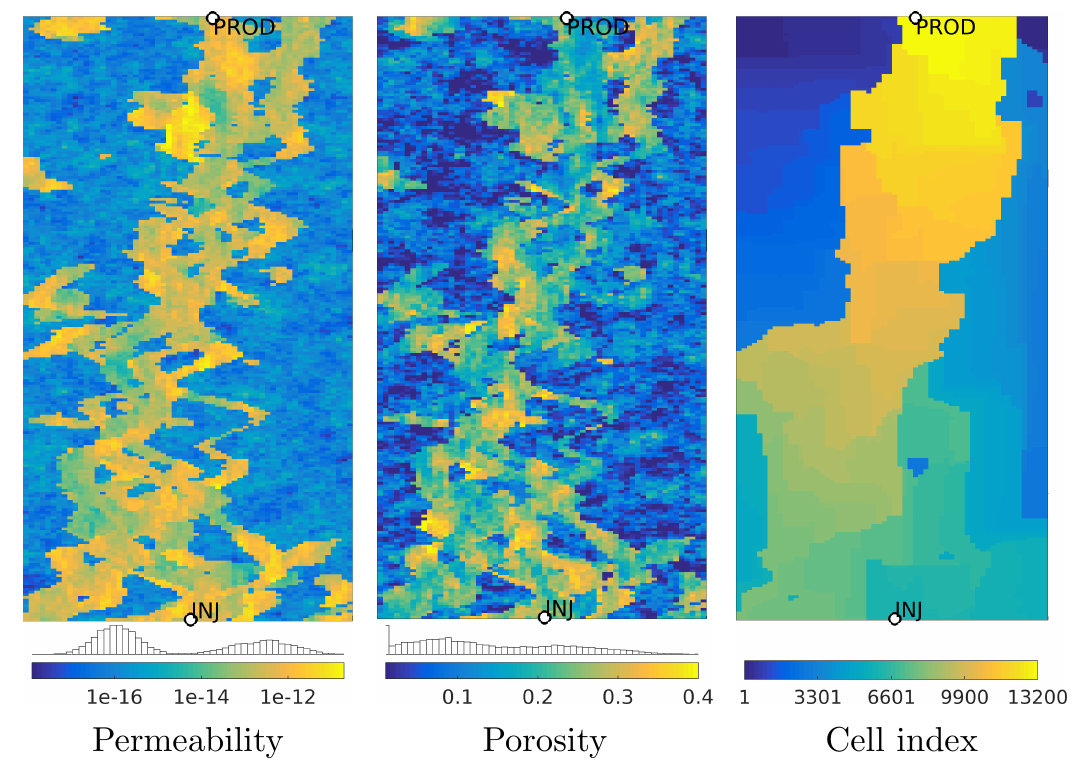}
    \caption{Petrophysical properties for Layer 50 from the SPE10 Benchmark. The right plot shows the cell index after reordering.}
    \label{fig:spe10-setup}
\end{figure}

\begin{figure*}[p!]
    \centering
    \includegraphics[width=\textwidth]{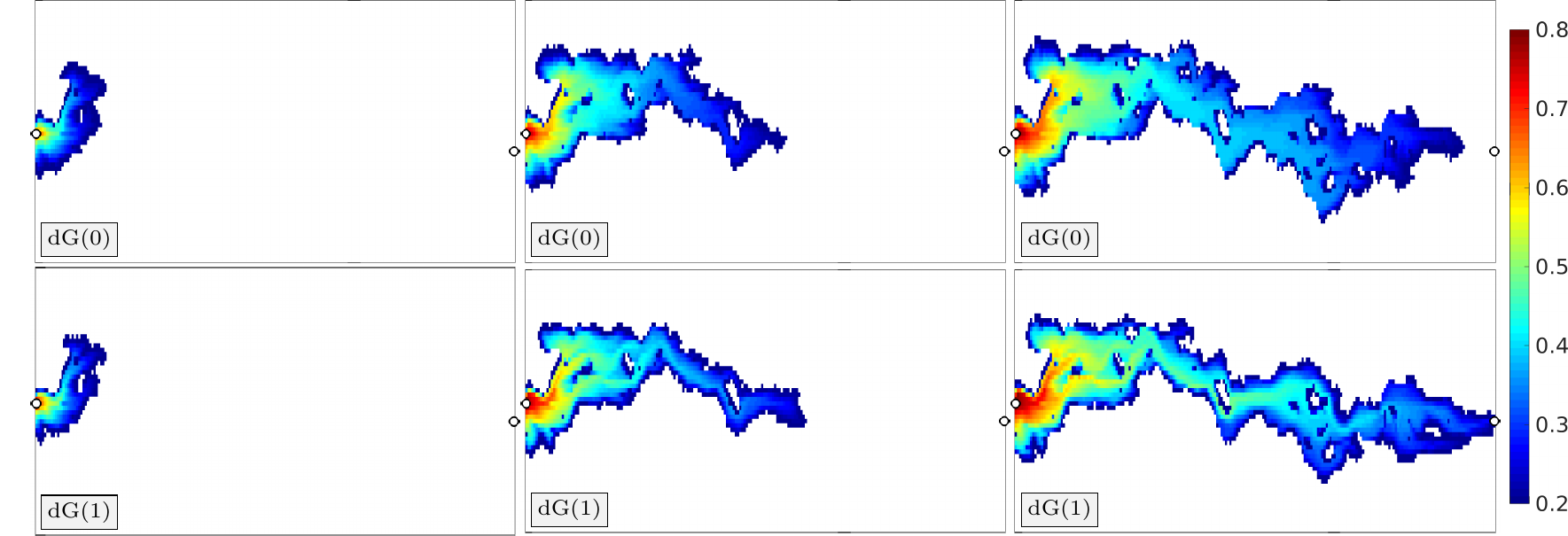}
    \caption{Evolution of the water saturation for three different times in the channelized reservoir in Example~\ref{sec:spe10}. The upper row is for dG(0) and the lower row for dG(1).}
    \label{fig:spe10-sat}
    \vskip\floatsep
    \centering
    \includegraphics[width=\textwidth]{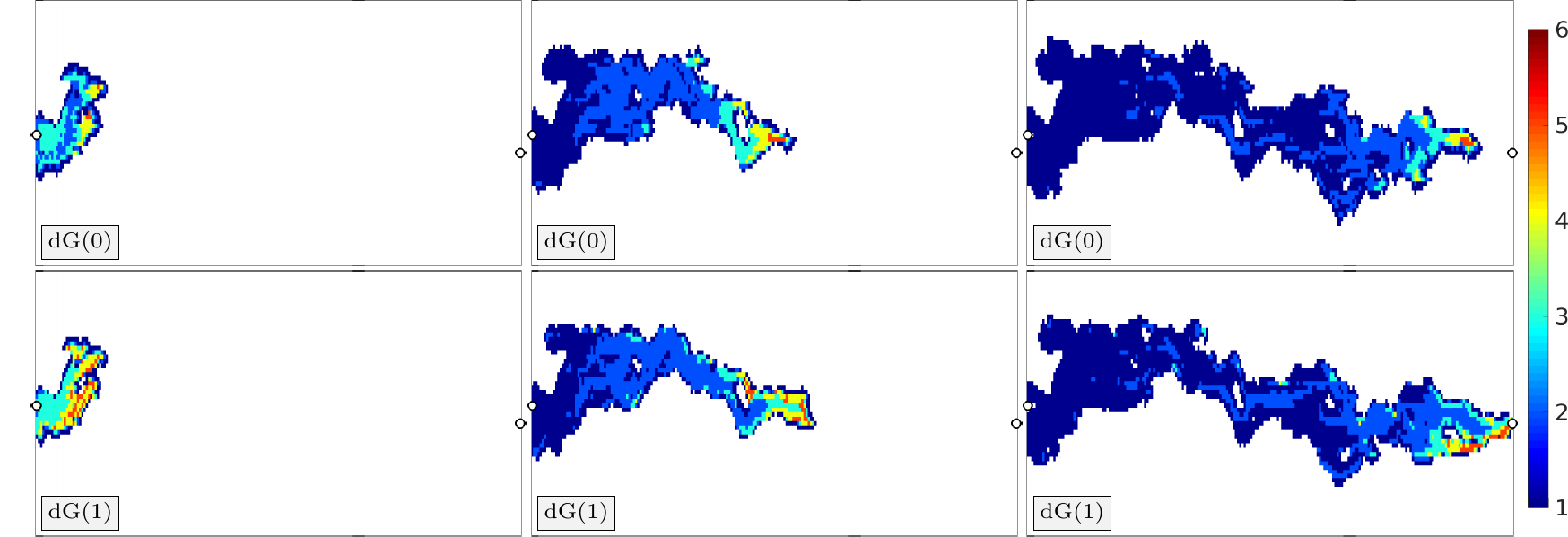}
    \caption{Number of nonlinear iterations per cell at three instances in time corresponding to the solution profiles shown in Figure~\ref{fig:spe10-sat}. Zero iterations were performed for cells in white. The upper row is for dG(0) and the lower row for dG(1).}
    \label{fig:spe10-it}
    \vskip\floatsep    
    \hspace{-0em}\includegraphics[width=0.95\textwidth]{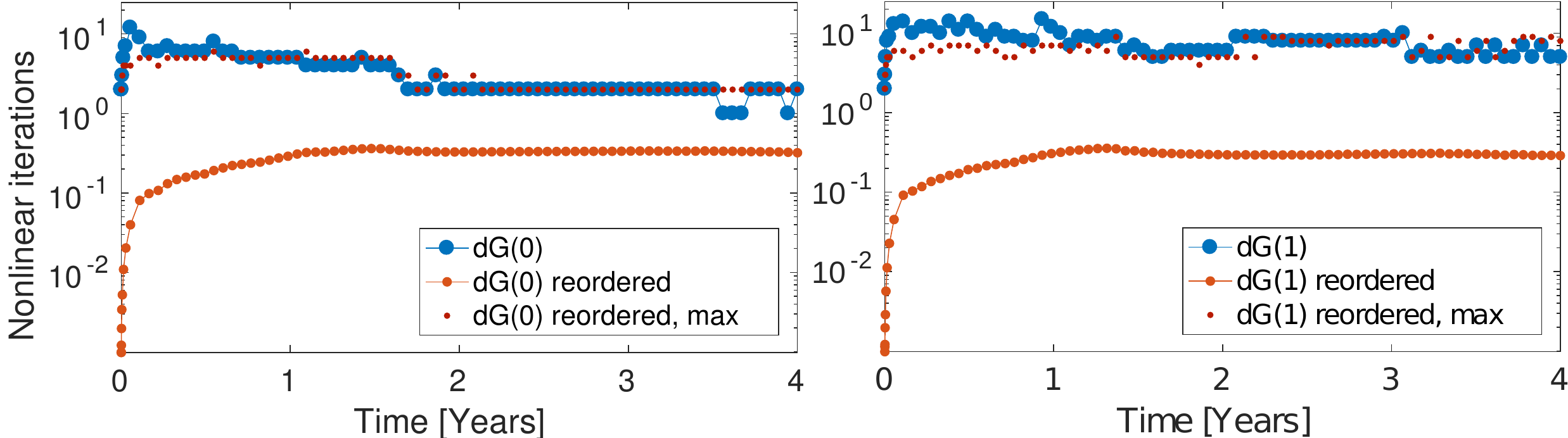}    
    \caption{Comparison of average number of nonlinear iterations per cell as function of time for Example~\ref{sec:spe10}. Note that the $y$-axis is logarithmic.}
    \label{fig:spe10-it-graph}
\end{figure*}

\begin{figure*}[t!]
    \centering
    \includegraphics[width=1\textwidth]{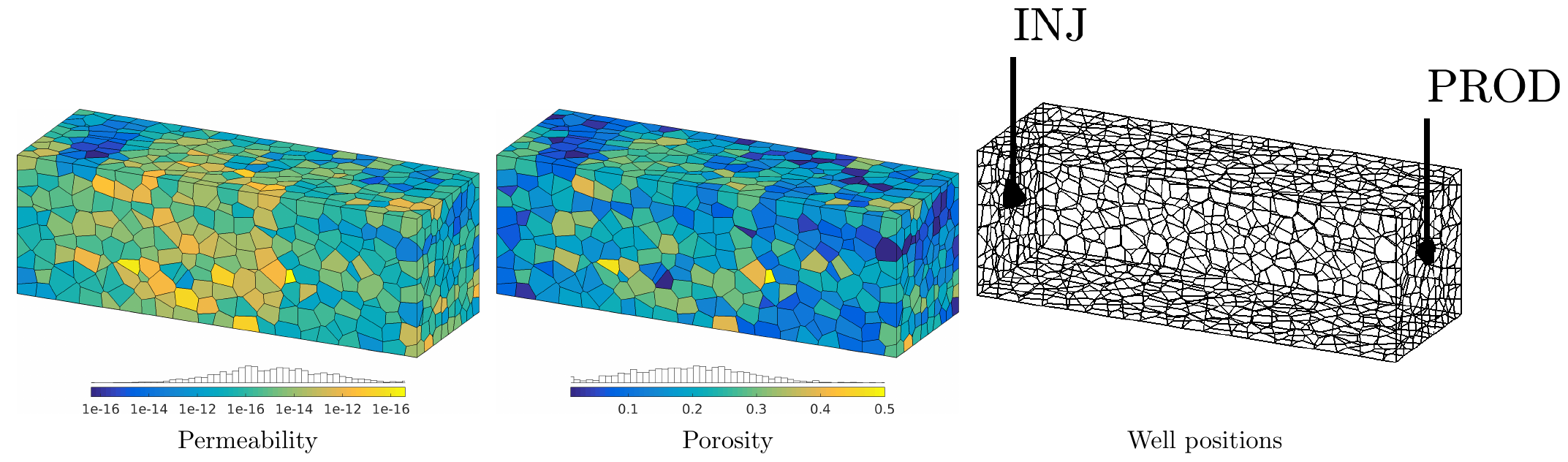}
    \caption{Petrophysical properties and well setup for the Voronoi-example.}
    \label{fig:pebi-setup}
\end{figure*}

Next, we consider a 2D simulation model with petrophysical properties sampled from Layer 50 from the SPE10 Model 2 Benchmark \citep{SPE10}. This layer is part of a fluvial formation and consists of an intertwined pattern of high-permeable sand channels on a background of low-permeable mudstone. We inject 0.2 pore volumes of water at a constant rate in one end of a channel, and produce at a fixed bottom-hole pressure of 275 bars at the opposite end of the same channel. The reservoir is initially filled with a mixture of water and oil, and the rock and fluids are weakly compressible. Water and oil have quadratic relative permeabilities with residual saturations of 0.2, and viscosities of 2.85 and 3.0 cP, respectively. Permeability and porosity are depicted in Figure~\ref{fig:spe10-setup}. The figure also shows the cell index after reordering. As expected, the cell number increases from the injector and along the high-flow channels to the producer. Interestingly, the model contains cycles of cells that are \textit{not} connected to the injector--producer pair by nonzero fluxes. This is an effect of fluid compressibility, and the reordering algorithm assigns lower indices to cells which essentially should remain inactive throughout the transport time step. 

We simulate using dG(0) and dG(1) with reordering implemented in MRST. Figure~\ref{fig:spe10-sat} shows the resulting water saturation at selected time steps, whereas Figure~\ref{fig:spe10-it} shows the iterations used per active cell in the corresponding time steps. By active cells, we mean cells in which the water saturation was updated during the time step either because the residual at the beginning of the time step was above the nonlinear tolerance, or because one or more of the cell's upstream neighbors were updated during the time step. Notice that iterations are performed only in a small fraction of the cells in each time step and that the updates are chiefly located to the fluvial channel. We also see that the higher-order solution profile is much less diffuse. Effectively, the water channel fills faster, so that dG(1) predicts earlier water breakthrough than dG(0).

For comparison, Figure~\ref{fig:spe10-it-graph} also reports the (average) number of nonlinear Newton--Raphson iterations used to solve each time step. For the reordered solver, we report the average number of nonlinear iterations used per cell, defined over all cells, and the maximum number of iterations observed in the cell with slowest nonlinear convergence. Because large fractions of the cells remain inactive in most of the time steps, the average number of cellwise nonlinear iterations is far below one in all time steps. This confirms previous observations by \citet{Natvig2008}. Interestingly, we see that the maximum number of nonlinear iterations for the reordered solver in most time steps is \textit{less than} the number of Newton--Raphson solves used when solving for all cells simultaneously for the same time step. This indicates that the iteration path followed by individual cells in a Newton--Raphson solve of all cells simultaneously is not necessarily optimal.

\begin{figure*}[t!]
    \centering
    \includegraphics[width=\textwidth]{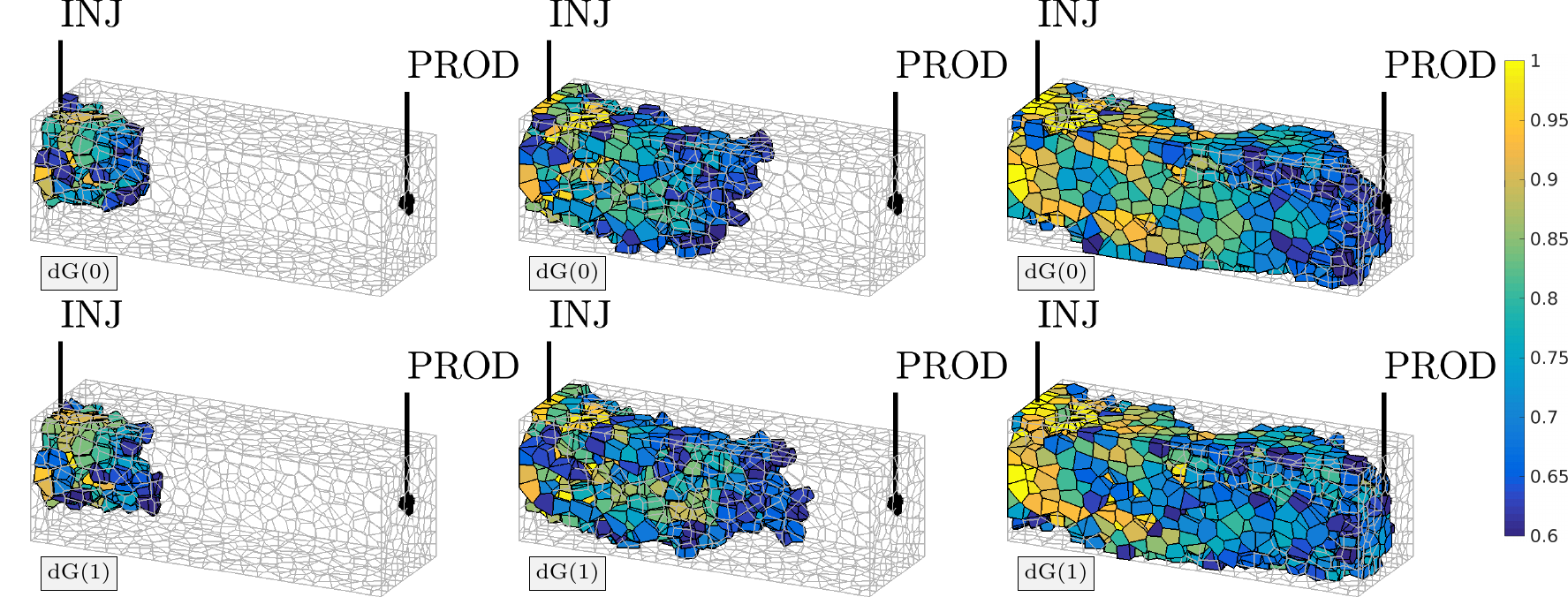}
    \caption{Thresholded water saturation profiles for the two dG solutions at three selected time steps for Example~\ref{sec:voronoi}.}
    \label{fig:pebi-sat}
    \vskip\floatsep
    \includegraphics[width = 0.95\textwidth]{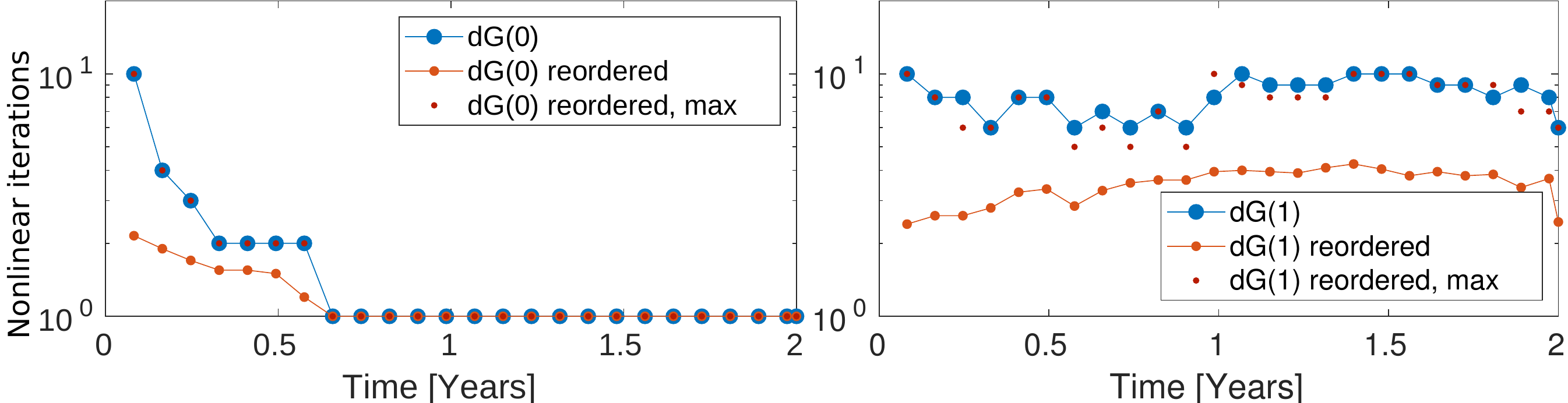}
    \caption{Average number of cellwise nonlinear iterations per time step used for each of the solvers in Example~\ref{sec:voronoi}, with dG(0) in the left plot and dG(1) in the right plot.}
    \label{fig:pebi-its}
\end{figure*}

\subsection{3D Voronoi grid}
\label{sec:voronoi}
In this example, we study the effect of using higher-order dG methods and reordering to simulate a case posed on a 3D, fully unstructured PEBI/Voronoi grid generated using the \texttt{upr} module in MRST \citep{Berge2018, rsc_grid_disc}. The domain is comprised of a rectangular box of dimensions $300 \times 100 \times 100$ m$^3$ with approximately 20, 10, and 10 cells in each of the axial directions, respectively. Porosity and permeability are sampled from the top layers of Model 2 from SPE10 \citep{SPE10}. The injector is placed in a cell at the center of the left boundary. We inject water at constant bottom-hole pressure of 600 bars. A producer operated at a fixed bottom-hole pressure of 50 bars is placed at the center of the right boundary. Figure~\ref{fig:pebi-setup} shows petrophysical properties and well positions.

\begin{figure}[h!]
    \centering
    \includegraphics[width = 0.46\textwidth]{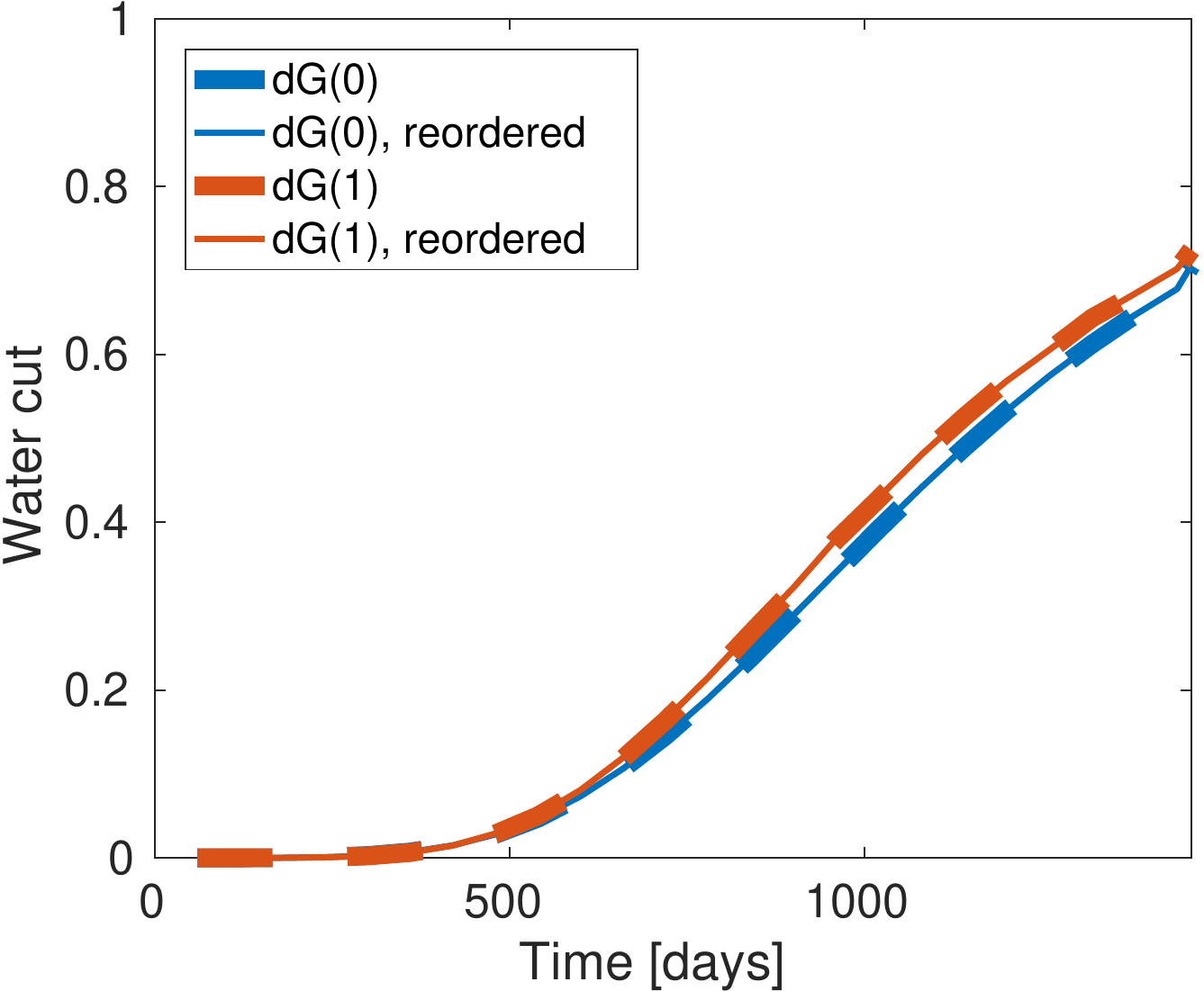}
    \caption{Water cut in the production well in Example~\ref{sec:voronoi}.}
    \label{fig:pebi-wcut}
\end{figure}

We simulate the case using the dG(0) and dG(1) solvers from MRST with and without reordering. For the reordering algorithm, we use the block version with a block size $n_b$ of 100 cells. Figure~\ref{fig:pebi-sat} reports saturation profiles at three selected time steps. Compared with dG(1), the dG(0) saturation profiles are more diffusive and fill more of the reservoir cross-section. As a result, dG(1) predicts a higher water-cut in the producer, as reported in Figure~\ref{fig:pebi-wcut}.

The solutions computed with the ordered and the original Newton--Raphson solvers are identical to plotting resolution, which verifies the correctness of the localization and the new block algorithm. Finally, we look at the number of iterations used by the different solvers, shown in Figure~\ref{fig:pebi-its}. We observe that even with a block size of $n_b = 100$ cells, the average number of nonlinear iterations per cell is significantly smaller than the number of iterations used when solving for all cells simultaneously.

\section{Concluding remarks}

The transport substep in a sequentially implicit solution procedure for black-oil models can be significantly accelerated by ordering the grid cells in the simulation model according to total velocity computed in a proceeding pressure step, so that the residual transport equations can be solved cell-by-cell from inflow to outflow. This not only localizes the nonlinear Newton solver, so that the computational effort can be focused on cells where the transported quantities change significantly during each time step, but also contributes to improve the general robustness of the Newton solver.  As a verification of this concept in an industry-grade setting, we have used a sequential simulator from the open-source OPM initiative to demonstrate that a reordered, global nonlinear Gauss--Seidel strategy computes an acceptable solution for the full simulation model of the Norne oil field. Additional verifications have also been obtained using similar industry-grade simulators from MRST.

The ordering approach is also applicable to transport solvers based on discontinuous Galerkin discretizations in combination with single-point upstream mobility weighting. Such discretizations are detailed herein for general 3D cases with (almost) no assumptions on the polyhedral cell geometries and grid topology. By reordering cells optimally, one can ensure that the added computational cost of solving for more degrees of freedom to obtain higher spatial order remains local to each cell. Our results obtained using MRST show that a second-order dG(1) method has significantly improved accuracy compared with the standard first-order volume method used in industry at the expense of only a moderate increase in the number of nonlinear iterations.

Herein, we have only presented simulations of two-phase oil-water and three-phase black-oil systems, but the same principles are applicable to a wide variety of different black-oil type and compositional models used to described secondary and tertiary recovery processes. As an example, \citet{klemetsdal2019:rsc} discuss several cases with significantly more complex flow physics, including compositional flow and gas injection with strong gravity effects, as well as use of dynamically adaptive grid coarsening.

\section{Acknowledgements}
The work of Klemetsdal, Rasmussen, and Lie is funded by the Research Council of Norway through grant no. 244361. M{\o}yner is funded by VISTA, a basic research program funded by Equinor and conducted in close collaboration with the Norwegian Academy of Science and Letters.

\begin{small}
  \bibliography{paper}

\begin{thebibliography}{34}
\providecommand{\natexlab}[1]{#1}
\providecommand{\url}[1]{\texttt{#1}}
\expandafter\ifx\csname urlstyle\endcsname\relax
  \providecommand{\doi}[1]{doi: #1}\else
  \providecommand{\doi}{doi: \begingroup \urlstyle{rm}\Url}\fi

\bibitem[Appleyard and Cheshire(1982)]{Appleyard1982}
J.~R. Appleyard and I.~M. Cheshire.
\newblock The cascade method for accelerated convergence in implicit
  simulators.
\newblock In \emph{European Petroleum Conference, 25-28 October, London, United
  Kingdom}. Society of Petroleum Engineers, October 1982.
\newblock \doi{10.2118/12804-MS}.

\bibitem[Bell et~al.(1986)Bell, Trangenstein, and Shubin]{Bell1986}
J.~B. Bell, J.~A. Trangenstein, and G.~R. Shubin.
\newblock Conservation laws of mixed type describing three-phase flow in porous
  media.
\newblock \emph{SIAM Journ. Appl. Math.}, 46\penalty0 (6):\penalty0 1000--1017,
  1986.
\newblock \doi{10.1137/0146059}.

\bibitem[Berge et~al.(2019)Berge, Klemetsdal, and Lie]{Berge2018}
R.~L. Berge, {\O}.~S. Klemetsdal, and K.-A. Lie.
\newblock Unstructured voronoi grids conforming to lower dimensional objects.
\newblock \emph{Comput. Geosci.}, 23\penalty0 (1):\penalty0 169--188, 2019.
\newblock \doi{10.1007/s10596-018-9790-0}.

\bibitem[Bratvedt et~al.(1996)Bratvedt, Gimse, and Tegnander]{Bratvedt1996}
F.~Bratvedt, T.~Gimse, and C.~Tegnander.
\newblock Streamline computations for porous media flow including gravity.
\newblock \emph{Transp. Porous Media}, 25\penalty0 (1):\penalty0 63--78, 1996.
\newblock \doi{10.1007/BF00141262}.

\bibitem[Brenier and Jaffr{\'e}(1991)]{brenier1991upstream}
Y.~Brenier and J.~Jaffr{\'e}.
\newblock Upstream differencing for multiphase flow in reservoir simulation.
\newblock \emph{SIAM J. Numer. Anal.}, 28\penalty0 (3):\penalty0 685--696,
  1991.
\newblock \doi{10.1137/0728036}.

\bibitem[Christie and Blunt(2001)]{SPE10}
M.~A. Christie and M.~J. Blunt.
\newblock Tenth {SPE} comparative solution project: A comparison of upscaling
  techniques.
\newblock \emph{SPE Reservoir Eval.~Eng.}, 4:\penalty0 308--317, 2001.
\newblock \doi{10.2118/72469-PA}.

\bibitem[Cockburn and Shu(1989)]{Cockburn1989}
B.~Cockburn and C.-W. Shu.
\newblock {TVB} runge-kutta local projection discontinuous galerkin finite
  element method for conservation laws ii: General.
\newblock \emph{Math. Comp.}, 52\penalty0 (186):\penalty0 411--435, 1989.
\newblock \doi{10.1090/S0025-5718-1989-0983311-4}.

\bibitem[Cockburn and Shu(1991)]{Cockburn1991}
B.~Cockburn and C.-W. Shu.
\newblock The runge-kutta local projection $p^1$- discontinuous-galerkin finite
  element method for scalar conservation laws.
\newblock \emph{ESAIM--Math. Model. Num.}, 25\penalty0 (3):\penalty0 337--361,
  1991.

\bibitem[Datta-Gupta and King(2007)]{Datta-Gupta2007}
A.~Datta-Gupta and M.~J. King.
\newblock \emph{Streamline simulation: theory and practice}, volume~11 of
  \emph{SPE Textbook Series}.
\newblock Society of Petroleum Engineers, 2007.

\bibitem[Eikemo et~al.(2009)Eikemo, Lie, Dahle, and Eigestad]{Eikemo09}
B.~Eikemo, K.-A. Lie, H.~K. Dahle, and G.~T. Eigestad.
\newblock Discontinuous {Galerkin} methods for transport in advective transport
  in single-continuum models of fractured media.
\newblock \emph{Adv. Water Resour.}, 32\penalty0 (4):\penalty0 493--506, 2009.
\newblock \doi{10.1016/j.advwatres.2008.12.010}.

\bibitem[Gries et~al.(2014)Gries, St{\"u}ben, Brown, Chen, and
  Collins]{Gries2014}
S.~Gries, K.~St{\"u}ben, G.~L. Brown, D.~Chen, and D.~A. Collins.
\newblock Preconditioning for efficiently applying algebraic multigrid in fully
  implicit reservoir simulations.
\newblock \emph{SPE J.}, 19\penalty0 (04):\penalty0 726--736, 2014.
\newblock \doi{10.2118/163608-PA}.

\bibitem[Jenny et~al.(2006)Jenny, Lee, and Tchelepi]{Jenny2006}
P.~Jenny, S.~H. Lee, and H.~A. Tchelepi.
\newblock Adaptive fully implicit multi-scale finite-volume method for
  multi-phase flow and transport in heterogeneous porous media.
\newblock \emph{J. Comput. Phys.}, 217\penalty0 (2):\penalty0 627--641, 2006.
\newblock \doi{10.1016/j.jcp.2006.01.028}.

\bibitem[Klausen et~al.(2012)Klausen, Rasmussen, and Stephansen]{Klausen2012}
R.~A. Klausen, A.~F. Rasmussen, and A.~Stephansen.
\newblock Velocity interpolation and streamline tracing on irregular
  geometries.
\newblock \emph{Comput. Geosci.}, 16:\penalty0 261--276, 2012.
\newblock \doi{10.1007/s10596-011-9256-0}.

\bibitem[Klemetsdal et~al.(2017)Klemetsdal, Berge, Lie, Nilsen, and
  M{\o}yner]{rsc_grid_disc}
{\O}.~S. Klemetsdal, R.~L. Berge, K.-A. Lie, H.~M. Nilsen, and O.~M{\o}yner.
\newblock Unstructured gridding and consistent discretizations for reservoirs
  with faults and complex wells.
\newblock In \emph{SPE Reservoir Simulation Conference}. Society of Petroleum
  Engineers, 2017.

\bibitem[Klemetsdal et~al.(2019{\natexlab{a}})Klemetsdal, M{\o}yner, and
  Lie]{Klemetsdal2018}
{\O}.~S. Klemetsdal, O.~M{\o}yner, and K.-A. Lie.
\newblock Robust nonlinear newton solver with adaptive interface-localized
  trust regions.
\newblock \emph{SPE J.}, 2019{\natexlab{a}}.
\newblock \doi{10.2118/195682-PA}.

\bibitem[Klemetsdal et~al.(2019{\natexlab{b}})Klemetsdal, M{\o}yner, and
  Lie]{klemetsdal2019:rsc}
{\O}.~S. Klemetsdal, O.~M{\o}yner, and K.-A. Lie.
\newblock Implicit high-resolution compositional simulation with optimal
  ordering of unknowns and adaptive spatial refinement.
\newblock In \emph{SPE Reservoir Simulation Conference, 10-11 April, Galveston,
  Texas, USA}. Society of Petroleum Engineers, 2019{\natexlab{b}}.
\newblock \doi{10.2118/193934-MS}.

\bibitem[Kwok and Tchelepi(2007)]{Kwok2007}
F.~Kwok and H.~Tchelepi.
\newblock Potential-based reduced newton algorithm for nonlinear multiphase
  flow in porous media.
\newblock \emph{J. Comput. Phys.}, 227\penalty0 (1):\penalty0 706 -- 727, 2007.
\newblock \doi{10.1016/j.jcp.2007.08.012}.

\bibitem[Lie(2016)]{Lie2016MRST-book}
K.-A. Lie.
\newblock \emph{An Introduction to Reservoir Simulation Using {MATLAB/GNU
  Octave}: User guide for the {Matlab Reservoir Simulation Toolbox (MRST)}}.
\newblock Cambridge University Press, 2016.

\bibitem[Lie et~al.(2012{\natexlab{a}})Lie, Krogstad, Ligaarden, Natvig,
  Nilsen, and Skaflestad]{MRST12:comg}
K.-A. Lie, S.~Krogstad, I.~S. Ligaarden, J.~R. Natvig, H.~M. Nilsen, and
  B.~Skaflestad.
\newblock Open source {MATLAB} implementation of consistent discretisations on
  complex grids.
\newblock \emph{Comput. Geosci.}, 16:\penalty0 297--322, 2012{\natexlab{a}}.
\newblock ISSN 1420-0597.
\newblock \doi{10.1007/s10596-011-9244-4}.

\bibitem[Lie et~al.(2012{\natexlab{b}})Lie, Natvig, and Nilsen]{Lie2012}
K.~A. Lie, J.~R. Natvig, and H.~M. Nilsen.
\newblock Discussion of dynamics and operator splitting techniques for
  two-phase flow with gravity.
\newblock \emph{Int. J. Numer. Anal. Mod.}, 9\penalty0 (3):\penalty0 684--700,
  2012{\natexlab{b}}.

\bibitem[Lie et~al.(2014)Lie, Nilsen, Rasmussen, and Raynaud]{Lie2014}
K.-A. Lie, H.~M. Nilsen, A.~F. Rasmussen, and X.~Raynaud.
\newblock Fast simulation of polymer injection in heavy-oil reservoirs on the
  basis of topological sorting and sequential splitting.
\newblock \emph{SPE J.}, 19\penalty0 (06):\penalty0 0991--1004, 2014.
\newblock \doi{10.2118/163599-PA}.

\bibitem[Lie et~al.(2017)Lie, M{\o}yner, Natvig, Kozlova, Bratvedt, Watanabe,
  and Li]{Lie2017ms-review}
K.-A. Lie, O.~M{\o}yner, J.~R. Natvig, A.~Kozlova, K.~Bratvedt, S.~Watanabe,
  and Z.~Li.
\newblock Successful application of multiscale methods in a real reservoir
  simulator environment.
\newblock \emph{Comput. Geosci.}, 21\penalty0 (5):\penalty0 981--998, Dec 2017.
\newblock \doi{10.1007/s10596-017-9627-2}.

\bibitem[M{\o}yner(2017)]{Moyner2016}
O.~M{\o}yner.
\newblock Nonlinear solver for three-phase transport problems based on
  approximate trust regions.
\newblock \emph{Comput. Geosci.}, 21\penalty0 (5-6):\penalty0 999--1021, 2017.
\newblock \doi{10.1007/s10596-017-9660-1}.

\bibitem[M{\"u}ller et~al.(2013)M{\"u}ller, Kummer, and Oberlack]{Muller2013}
B.~M{\"u}ller, F.~Kummer, and M.~Oberlack.
\newblock Highly accurate surface and volume integration on implicit domains by
  means of moment-fitting.
\newblock \emph{Int. J. Numer. Meth. Eng.}, 96\penalty0 (8):\penalty0 512--528,
  2013.
\newblock \doi{10.1002/nme.4569}.

\bibitem[Natvig and Lie(2008)]{Natvig2008}
J.~R. Natvig and K.-A. Lie.
\newblock Fast computation of multiphase flow in porous media by implicit
  discontinuous {Galerkin} schemes with optimal ordering of elements.
\newblock \emph{J. Comput. Phys.}, 227\penalty0 (24):\penalty0 10108--10124,
  2008.
\newblock \doi{10.1016/j.jcp.2008.08.024}.

\bibitem[Natvig et~al.(2007)Natvig, Lie, Eikemo, and Berre]{Natvig07}
J.~R. Natvig, K.-A. Lie, B.~Eikemo, and I.~Berre.
\newblock An efficient discontinuous {Galerkin} method for advective transport
  in porous media.
\newblock \emph{Adv. Water Resour.}, 30\penalty0 (12):\penalty0 2424--2438,
  2007.
\newblock \doi{10.1016/j.advwatres.2007.05.015}.

\bibitem[OPM(2019)]{opm}
OPM.
\newblock The open porous media {(OPM)} initative.
\newblock \url{https://opm-project.org/}, 2019.

\bibitem[Rasmussen and Lie(2014)]{Rasmussen2014}
A.~F. Rasmussen and K.-A. Lie.
\newblock Discretization of flow diagnostics on stratigraphic and unstructured
  grids.
\newblock In \emph{16th European Conference on the Mathematics of Oil Recovery,
  Catalania, Sicily, Italy}. European Association of Geoscientits and
  Engineers, 2014.

\bibitem[Shahvali and Tchelepi(2013)]{Shahvali2013}
M.~Shahvali and H.~A. Tchelepi.
\newblock Efficient coupling for non-linear multiphase flow with strong
  gravity.
\newblock In \emph{SPE Reservoir Simulation Symposium, 18-20 February, The
  Woodlands, Texas, USA}. Society of Petroleum Engineers, February 2013.
\newblock \doi{10.2118/163659-MS}.

\bibitem[Sheth and Younis(2017)]{Sheth2017}
S.~M. Sheth and R.~M. Younis.
\newblock Localized solvers for general full-resolution implicit reservoir
  simulation.
\newblock In \emph{SPE Reservoir Simulation Conference}. Society of Petroleum
  Engineers, 2017.
\newblock \doi{10.2118/182691-MS}.

\bibitem[{The Open Porous Media (OPM) Initative}(2015)]{norne-data}
{The Open Porous Media (OPM) Initative}.
\newblock The {Norne} dataset.
\newblock \url{https://github.com/OPM/opm-data}, 2015.

\bibitem[Trangenstein and Bell(1989)]{Trangenstein1989}
J.~A. Trangenstein and J.~B. Bell.
\newblock Mathematical structure of the black-oil model for petroleum reservoir
  simulation.
\newblock \emph{SIAM J. Appl. Math.}, 49\penalty0 (3):\penalty0 749--783, 1989.
\newblock \doi{10.1137/0149044}.

\bibitem[Trottenberg et~al.(2000)Trottenberg, Oosterlee, and
  Schuller]{Trottenberg2000}
U.~Trottenberg, C.~W. Oosterlee, and A.~Schuller.
\newblock \emph{Multigrid}.
\newblock Academic press, 2000.

\bibitem[Watts(1986)]{Watts1986}
J.~Watts.
\newblock A compositional formulation of the pressure and saturation equations.
\newblock \emph{SPE Reservoir Eng.}, 1\penalty0 (03):\penalty0 243--252, 1986.
\newblock \doi{10.2118/12244-PA}.

\end{thebibliography}
\end{small}

\end{document}